\documentclass[final,1p,times]{elsarticle}

\usepackage{amssymb}
\usepackage{graphicx}
\usepackage{amsthm}
\usepackage{amsfonts,amsmath}
\usepackage{subcaption}
\usepackage{cleveref}
\usepackage{url}
\usepackage{float} 
\usepackage{color,changebar}
\usepackage{algorithm}
\usepackage{multirow}
\usepackage[noend]{algpseudocode}
\usepackage[utf8]{inputenc}
\usepackage{circuitikz}
\usepackage{comment}
\usepackage{chemmacros}
\theoremstyle{remark}
\newtheorem{remark}{Remark}

\begin{document}

\begin{frontmatter}


\title{Parametric Dynamic Mode Decomposition for nonlinear parametric dynamical systems}

\author[a]{Shuwen Sun\corref{cor1}}
\ead{ssun@mpi-magdeburg.mpg.de}
\author[a]{Lihong Feng}
\author[b]{Hoon Seng Chan}
\author[a]{Tamara Mili\v{c}i{\'c}}
\author[a]{Tanja Vidakovi{\'c}-Koch}
\author[c]{Fridolin R{\"o}der}
\author[a,d]{Peter Benner}

\affiliation[a]{organization={Max Planck Institute for Dynamics of Complex Technical Systems},
             city={Magdeburg},
             postcode={39016},
             country={Germany}}
             
\affiliation[b]{organization={Karlsruhe Institute of Technology, Institute for Applied Materials –
Electrochemical Technologies},
             city={Karlsruhe},
             postcode={76131},
             country={Germany}}
             
\affiliation[c]{organization={Bavarian Center for Battery Technology (BayBatt), University of Bayreuth},
			city={Bayreuth}, 
			postcode={95447},
			country={Germany}}
             
\affiliation[d]{organization={Faculty of Mathematics, Otto von Guericke University},
             city={Magdeburg},
             postcode={39016},
             country={Germany}}

\cortext[cor1]{Corresponding author at Max Planck Institute for Dynamics of Complex Technical Systems, Sandtorstraße 1, 39106, Magdeburg, Germany.}


%
%

\begin{abstract}

A non-intrusive model order reduction (MOR) method that combines features of the dynamic mode decomposition (DMD) and the radial basis function (RBF) network is proposed to predict the dynamics of parametric nonlinear systems. In many applications, we have limited access to the information of the whole system, which motivates non-intrusive model reduction. One bottleneck is capturing the dynamics of the solution without knowing the physics inside the ``black-box'' system. DMD is a powerful tool to mimic the dynamics of the system and to give a reliable approximation of the solution in the time domain using only the dominant DMD modes. However, DMD in general cannot reproduce the parametric behavior of the dynamics.
Our contribution focuses on extending DMD to parametric DMD by RBF interpolation. Specifically, a RBF network is first trained using snapshot matrices at a limited number of parameter samples. The snapshot matrix at any new parameter sample can be quickly learned from the RBF network. DMD then uses the newly generated snapshot matrix at the online stage to predict the time patterns of the dynamics corresponding to the new parameter sample. The proposed framework and algorithm are tested and validated by numerical examples including models with parametrized and time-varying inputs.

\end{abstract}

%


\end{frontmatter}


\section{Introduction}%
\label{sec:intro}
Nonlinear dynamical systems arise from many physical and engineering applications. Solving systems with nonlinear effects and parameter variations indeed costs a lot of time and effort, which motivates model order reduction, a technique for constructing compact surrogates of nonlinear systems to realize accelerated computation with acceptable accuracy. The computational efforts in constructing the surrogate, i.e., the reduced-order model (ROM), is usually concentrated at the offline stage, while the process of employing the ROM for simulation or any other multi-query tasks is known as the online stage. When the online stage is fast enough, it can be stated as ``real-time'' computation and is promising for real applications. There are various subtopics and methods in MOR aiming at different applications, such as modal truncation, balanced truncation \cite{morGugA04, morMehS05}, Krylov subspace methods (moment matching) \cite{morFreR03}, local linear embedding (LLE) \cite{morRowS00}, proper orthogonal decomposition (POD, also known as principal component analysis in the statistical area or Karhunen-Loeve expansion in the stochastic area) \cite{Pea01, Lum67, Lum81, morSir87}, reduced basis methods, dynamic mode decomposition (DMD) \cite{morTuRLetal14}, data-driven and machine learning approaches \cite{morZanGT21, morGusS99, morGriG15,  morNakST18, morXiaGTetal19, morRegDQ20, morRahPSetal19, morBhaHKetal21, morRenMR20, morKasGH20,  morBenGKetal20, morGoyB21}.

When a dynamic system is seen as a ``black box" so that the only information of the system are the inputs and its corresponding outputs, intrusive MOR based on projection is impossible, and non-intrusive MOR is preferred. Efficient MOR for nonlinear time-evolution systems parametrized with some physical or geometrical parameters is challenging. Although intrusive MOR based on projection for such systems has achieved a great amount of success~\cite{morhandbookV1, morhandbookV2, morhandbookV3}, non-intrusive MOR methods that are robust for systems characterized by all the above three properties, i.e., nonlinear, parametric and time-dependent, are still not fully explored, though some methods are proposed~\cite{morBenGW15, morFreDM21, morRegDQ20, morGuoH19, morRenMR20, morKasGH20, morBenGKetal20, morCheWHetal20, morGoyB21, morXiaFN17, morXuD20}. At present, more and more non-intrusive MOR methods are based on machine learning to tackle such systems with strong nonlinearity~\cite{morFreDM21, morGuoH19, morRenMR20, morXiaFN17, morCheFetal20, morXuD20}. Furthermore, many of the existing methods assume that the solution space is of low dimension, and a global reduced space over the whole parameter domain is assumed~\cite{morRegDQ20, morRahPSetal19, morBhaHKetal21, morRenMR20, morGoyB21}. Fewer non-intrusive methods are successful for systems in which the solution is non-smooth in the parameter domain~\cite{morFreDM21, morXiaFN17, morCheFetal20, morXuD20}. Non-intrusive MOR methods with emphasis on treating non-smooth or convection-dominated problems are also proposed~\cite{morSarGB20, morMenBAetal20, morWel17, morHenMS22}. To the best of the authors' knowledge, many of them are only applicable to either parametric steady problems or time-evolution problems without parameters.

Dynamic mode decomposition can provide a way of discovering low-rank space-time patterns of the dynamics in an equation-free manner \cite{morBruK21}. DMD was first introduced to realize the nonlinear evolution of fluid dynamics. Based on the snapshot matrix from the system, DMD computes a linear operator that maps the snapshots one time step further. It appeared firstly in \cite{morSch10} and then it was later used for model order reduction. There exist different variants of the DMD method to overcome the different drawbacks of the standard DMD, such as Extended DMD \cite{morWilKR15} and Kernel DMD \cite{morWilRK15}. DMD is also combined with an autoencoder for non-intrusive model reduction of nonlinear dynamical systems \cite{OttR19}.

This work focuses on extending DMD to parametric DMD by combining DMD with the RBF network to achieve fast approximation of both the parametric behavior and time-evolution of the dynamics in a non-intrusive way. Compared with the existing methods based on deep learning, our proposed method is much faster to train, since the RBF network is known as a shallow neural network with much fewer parameters to be optimized during network training. Yet, the derived ROMs are still of acceptable accuracy. Some closely related methods are proposed in \cite{morXiaFN17}, where the RBF network is combined with POD and is also used for prediction in the time domain. Due to the limitation of the RBF interpolation only in the time domain, the method in~\cite{morXiaFN17} cannot predict the solution at a future time that is outside of the time interval used for training. Another recent work on parametric DMD~\cite{morHuhQetal2023} aims at reaching the same goal as our proposed method in different ways. Here, two different parametric DMD methods are proposed. The first method interpolates the eigenpairs associated with the projected Koopman matrices at different parameter samples. The second method instead interpolates the projected Koopman matrices corresponding to different parameter samples. Each method necessitates the second stage of interpolation: interpolating the associated DMD modes in the parameter domain to recover the solution in the original space. However, some limitations are also mentioned in this paper. The most restricting limitations are the following assumptions. Given the polynomial interpolation method used in \cite{morHuhQetal2023}, the smoothness of the eigenpairs over the parameter domain must be satisfied for the accuracy of the first method, and smoothness of the projected Koopman operator w.r.t the parameters is required for the second method to be successful. Another limitation lies in the fact that projected DMD used in \cite{morHuhQetal2023} can not assure that the dynamic modes are exactly the eigenvectors of the original Koopman matrix. In the latest paper \cite{morLuT22}, similar work has been done using DMD for non-intrusive MOR of parametric systems, the reduced-order model at any testing parameter sample is obtained from manifold-interpolation of the left singular vectors at training parameter samples and manifold-interpolation of the projected Koopman matrices. Some hyperparameters need to be heuristically tuned to achieve success, for example, the reference configuration $i_0$, which could lead to failure of the method if not optimally chosen. Furthermore, the proposed DMD method can only reconstruct the observables of the solution. The solution needs to be recovered by implementing an inverse mapping from the observables to the state space. For observables with a complex expression, it is unclear how the inverse mapping can be computed. 

In this work, the power of DMD for time-dependent problems is combined with the RBF network to derive a method that is robust for prediction in both the parameter domain and the time domain. When compared to the existing DMD-based methods for MOR of the parametric dynamical system, the RBF network that is applied for snapshot interpolation leads to the proposed parametric DMD method with much fewer constraints.

The remaining part of the work is organized as follows. In \Cref{sec:dmd}, a general overview of DMD is provided. The algorithm of the exact DMD and the kernel DMD are presented for use in the next sections. In \Cref{sec:work}, radial basis function (RBF) interpolation is shortly introduced. Then the proposed method, a practical algorithm, and some discussions are given. In \Cref{sec:examples}, three examples from real applications are presented to demonstrate the robustness of the proposed method. We conclude the work in \Cref{sec:conclu} with further outlook.

\section{Dynamic Mode Decomposition}%
\label{sec:dmd}
DMD is a non-intrusive MOR method for time-dependent systems. It provides a low-dimensional representation of the system solution via spatiotemporal decomposition of the dynamics. The main tool is the singular value decomposition (SVD) of a large data matrix and the eigendecomposition of a small projected data matrix. 
Suppose we have a nonlinear dynamic system of ordinary differential equations (ODEs):
\begin{equation}
	\label{eq:ode}
    \mathbf{u}^{\prime}(t) = g(\mathbf{u}(t)),
\end{equation}
where the state vector $u(t) \in \mathbb{R}^{n}$ , $g$: $ \mathbb{R}^{n} \to \mathbb{R}^{n}$ is a nonlinear operator. Applying an explicit time integration scheme to \cref{eq:ode} results in the following nonlinear evolution,

\begin{equation}
	\label{eq:ode2}
    \mathbf{u}_{i+1}=F(\mathbf{u}_i), \quad i=0,\ldots,m-1. 
\end{equation}
Note that $F$ may also depend on $\mathbf{u}_{i-1}$, etc. for a multi-step integration scheme. For simplicity of explanation, those dependencies are omitted here.

Consider the snapshot matrix $\mathbf{X}_0$ and the shifted snapshot matrix $\mathbf{X}_1$ as follows:

\begin{equation}
\label{eq:x0x1}
\mathbf{X}_0=\left[\begin{array}{cccc}
\mid & \mid & & \mid \\
\mathbf{u}_{0} & \mathbf{u}_{1} & \cdots & \mathbf{u}_{m-1} \\
\mid & \mid & & \mid
\end{array}\right] \in \mathbb{R}^{n\times m}, \\
\mathbf{X}_1=\left[\begin{array}{cccc}
\mid & \mid & & \mid \\
\mathbf{u}_{1} & \mathbf{u}_{2} & \cdots & \mathbf{u}_{m} \\
\mid & \mid& & \mid
\end{array}\right] \in \mathbb{R}^{n\times m}, 
\end{equation}
where $\mathbf{u}_i = \mathbf{u}(t_i)$, $i=0,...,m$, are state vectors at time $t_k$ within a certain time interval. They are also known as snapshots. DMD uses a linear time evolution to approximate the nonlinear evolution in \cref{eq:ode2}, i.e.

\begin{equation} 
	\label{eq:linop}
    \mathbf{X}_1 = \mathbf{K} \mathbf{X}_0.
\end{equation}

Then it finds the best fit $\mathbf{A}$ for the linear operator $\mathbf{K} \in \mathbb R^{n\times n}$. Mathematically, we have

\begin{equation}
\mathbf{A}=\underset{\mathbf{\tilde A \in \mathbb R^{n\times n}}}{\operatorname{argmin}}\left\|\mathbf{X}_1-\mathbf{\tilde A}\mathbf{X}_0\right\|_{F}=\mathbf{X}_1 \mathbf{X}_0^{\dagger},
\end{equation}
where $\left\|\cdot\right\|_{F}$ is the Frobenius norm and $^{\dagger}$ is the pseudo-inverse operator. When $\mathbf{X}_0$ and $\mathbf{X}_1$ are linearly consistent, i.e., whenever $\mathbf{X}_0 \mathbf{c}=0$, then $\mathbf{X}_1 \mathbf{c}=0$, then it is proved in \cite{morTuRLetal14} that $\mathbf{A}$ satisfies \cref{eq:linop}, i.e., $\mathbf{X}_1 = \mathbf{A} \mathbf{X}_0$. From the eigendecomposition of $\mathbf{A}$ we obtain the eigenvalues and eigenvectors of $\mathbf{A}$. The eigenvectors are also known as the DMD modes \cite{morTuRLetal14}. Reconstruction of the state can be done using these DMD modes and their evolution configured by the eigenvalues. Each eigenvalue represents the growth/decay rate (real part of the complex value) and oscillation with different frequencies (imaginary part of the value) of the corresponding mode. When $n$ is large, the eigendecomposition of $\mathbf{A}$ becomes inefficient. The practical algorithm of implementing DMD takes use of dimension reduction via SVD of the initial snapshot matrix $\mathbf{X}_0$ to compute the dominant DMD modes from the (truncated) left singular vectors $\mathbf{U	}$ and the eigendecomposition of the small projected matrix $\mathbf{\hat{A}} = \mathbf{U}^* \mathbf{A}\mathbf{U}$. \Cref{algorithm:exactdmd} presents the detailed procedure of the exact DMD algorithm first proposed in \cite{morTuRLetal14}.

The main difference between exact DMD and a previously proposed standard DMD (also known as projected DMD) lies in the way of computing the DMD modes. For the standard DMD, a DMD mode is computed from the matrix $\mathbf{U}$ of left singular vectors:
\begin{equation}
	\label{eq:proj_phi}
    \mathbf{\hat{\varphi}} = \mathbf{U} \mathbf{w},
\end{equation}
where $\mathbf{w}$ is an eigenvector of $\mathbf{\hat{A}}$, corresponding to an eigenvalue $\lambda$. However, for the exact DMD, the DMD mode $\varphi$ is defined as lying in the image of $\mathbf{X}_1$ instead of that of $\mathbf{X}_0$. It is computed as follows:
\begin{equation}
	\label{eq:exact_phi}
    \mathbf{\varphi} = \frac{1}{\lambda} \mathbf{X}_1 \mathbf{V} \mathbf{\Sigma}^{-1} \mathbf{w}.
\end{equation}

The aim of computing $\mathbf{\varphi}$ following \cref{eq:exact_phi} is to make sure that $\mathbf{\varphi}$ is the eigenvector of the original linear operator $\mathbf{A}$, i.e., $\mathbf{A} \mathbf{\varphi} = \lambda \mathbf{\varphi}$. This property is used in Step 6 in \Cref{algorithm:exactdmd} for the reconstruction of the dynamics.  Whereas, $\mathbf{\hat \varphi}$ in \cref{eq:proj_phi} doesn't meet such a requirement. A detailed explanation can be found in \cite{morTuRLetal14}.

After the DMD modes are computed in \Cref{algorithm:exactdmd}, the solution at any future time $t_i$ can be reconstructed from the DMD modes, and their initial amplitudes $\mathbf{b}$ computed based on the initial solution, see Steps 5-6 in	 \Cref{algorithm:exactdmd}.

\begin{remark}
The truncation in Step 2 of \Cref{algorithm:exactdmd} did not appear in the original exact DMD in \cite{morTuRLetal14} but was included in the exact DMD algorithm presented in \cite{morKutBBetal16} so that the computational cost of the eigendecomposition of $\mathbf{\hat{A}}$ is further reduced.  The truncation rank $r$ is determined according to the energy criteria: 
\begin{equation}
\label{eq:truncationtolerance}
\frac{\sum_{i=r+1}^{d}\sigma_i}{\sum_{i=1}^{d}\sigma_i} \leq \eta,
\end{equation} 
where $\eta$ is a tolerance decided by the user. This may introduce truncation errors, however, we found in the numerical tests that when $r \ll d$, the DMD still produces results with acceptable accuracy. Furthermore, once the truncation is introduced, the DMD modes computed in Step 5 are no longer the eigenvectors of $\mathbf{A}$.
\end{remark}

\begin{algorithm}[h]
    \begin{algorithmic}[1]
        \caption{Exact DMD \cite{morTuRLetal14, morKutBBetal16}} \label{algorithm:exactdmd}
        \State Collect the snapshots for the snapshot matrix pair $\{ \mathbf{X}_0, \mathbf{X}_1 \}$ in \cref{eq:x0x1}.
        \State Compute the compact SVD of the data snapshot matrix  $\mathbf{X}_0= \mathbf{U} \mathbf{\Sigma} \mathbf{V^T}, \mathbf{U} \in \mathbb R^{n\times d}, \mathbf{\Sigma} \in \mathbb R^{d\times d}, \mathbf{V} \in \mathbb R^{m \times d}$, $d \leq \min(m,n)$ is the rank of $\mathbf{X}_0$. Truncate and keep only the $r < d $ leading eigenvalues and the corresponding eigenvectors, so that $\mathbf{X}_0 \approx \mathbf{U}_r \mathbf{\Sigma}_r \mathbf{V}_r^{T}$, where $\mathbf{U}_r \in \mathbb{R}^{n \times r}, \mathbf{\Sigma}_r \in \mathbb{R}^{n \times r}$, and $\mathbf{V}_r \in \mathbb{R}^ {m\times r}$.
        \State Compute $\hat{\mathbf{A}}=\mathbf{U}_r^T\mathbf{A}\mathbf{U}_r $ by replacing $\hat{\mathbf{A}}$ with $\mathbf{X_1} \mathbf{X_0}^\dagger$, and $\mathbf{X_0}$ with its SVD in Step 2, i.e., $\hat{\mathbf{A}}=\mathbf{U}_r^{T} \mathbf{A} \mathbf{U}_r=\mathbf{U}_r^{T}\mathbf{X}_1 \mathbf{V}_r \Sigma^{-1} \mathbf{U}_r^{T} \mathbf{U}_r =\mathbf{U}_r^{T} \mathbf{X}_1 \mathbf{V}_r \mathbf{\Sigma}_r^{-1}$.
        \State Compute the eigendecomposition of $\hat{\mathbf{A}}$: $\hat{\mathbf{A}} \mathbf{W} = \mathbf{W} \mathbf{\Lambda}$, with $\mathbf{\Lambda} = diag(\lambda_1,\ldots,\lambda_r)$.
        \State Compute the DMD modes $\mathbf{\Phi}=[\mathbf{\varphi}_1, \ldots, \mathbf{\varphi}_r]$ by $\mathbf{\Phi}=\mathbf{X}_1 \mathbf{V}  \mathbf{\Sigma}_r^{-1} \mathbf{W}$. Given the initial solution $\mathbf{u}_0$ and suppose it can be represented by the DMD modes, i.e., $\mathbf{u}_0=\mathbf{\Phi} \mathbf{b}$, then the vector of coefficients  $\mathbf{b} =(b_1,\ldots,b_r)^T $ can be computed as $\mathbf{b}=\mathbf{\Phi}^{\dagger} \mathbf{u}_0$.
        \State Reconstruct the solution at any future time $t_i>0$ using the DMD modes: $\mathbf{u}_i=\mathbf{A}^i \mathbf{u}_0 = \sum_{k=1}^{r} \varphi_{k} b_k \lambda_k^i$.
    \end{algorithmic}
\end{algorithm}

\subsection{Extended and kernel DMD}
DMD uses a linear evolution scheme \cref{eq:linop} to approximate the nonlinear evolution \cref{eq:ode2}, which might cause big errors for some problems with strong nonlinearities. To improve the accuracy of DMD, extended DMD (EDMD) was proposed in \cite{morWilKR15}. Assuming that the state vector $u$ in \cref{eq:ode} can be spanned by $s$ eigenfunctions $\phi_k(u), k=1,\ldots,s$ of the Koopman operator $\mathbf{\mathcal{K}}$, i.e., 

\begin{equation}
\mathbf{u} = \sum\limits_{k=1}^s \mathbf{v}_k \phi_k(\mathbf{u}),
\end{equation}
then the nonlinear evolution \cref{eq:ode2} can be fully described by the Koopman operator via its eigenfunctions, eigenvalues and modes (see \cite{morWilKR15} for detailed derivation), i.e.,
\begin{equation}
\label{eq:exdmd_F2K}
\mathbf{F}(\mathbf{u})=\sum\limits_{k=1}^s \mathbf{v}_k (\mathbf{\mathcal{K}} \phi_k) (\mathbf{u}) =\sum\limits_{k=1}^s \lambda_k \mathbf{v}_k \phi_k (\mathbf{u}).
\end{equation}
Here, $\mathbf{\mathcal{K}}$ is the Koopman operator, $\phi_k(\mathbf{u})$ are the Koopman eigenfunctions, $\mathbf{v}_k$ are the Koopman modes, and $\lambda_k$ are the Koopman eigenvalues. Motivated by \cref{eq:exdmd_F2K}, EDMD tries to approximate the nonlinear evolution \cref{eq:ode2} via approximating the Koopman operator, its eigenfunctions and modes. The Koopman operator is approximated by using not only the data matrices $\mathbf{X}_0$, $\mathbf{X}_1$ but also a dictionary of functions of the state vector (observables) $\{\psi_1(\mathbf{u}), \psi_2(\mathbf{u}), \ldots, \psi_M(\mathbf{u})\}$, we can define a vector valued observable $\mathbf{\psi}(\mathbf{u}) = [\psi_1(\mathbf{u}) \quad \psi_2(\mathbf{u}) \quad \ldots \quad \psi_M(\mathbf{u})]$. Then the Koopman operator is supposed to be approximated by a finite dimensional matrix $\mathbf{\tilde K}$ with a residual term:
\begin{equation}
(\mathbf{\mathcal{K}} \theta)(\mathbf{u}) = \mathbf{\psi}(\mathbf{u})(\mathbf{\tilde K} \mathbf{a}) + r(\mathbf{u}),
\end{equation}
where $\mathbf{a}$ are the coefficients to construct a vector observable $\theta(\mathbf{u}) = \psi(\mathbf{u}) \mathbf{a}$  with a linear combination of $M$ components from $\mathbf{\psi}(\mathbf{u})$ and $r(\mathbf{u})$ is the residual term for the approximation. Full state observable can be obtained when $\theta(\mathbf{u}) = \mathbf{u}$.

To minimize this residual term, an objective function based on a single observable can be formulated as:
\begin{equation}
\mathbf{J} = \frac{1}{2}\sum\limits_{i=0}^{m-1} \lvert (\mathbf{\Psi}_1 - \mathbf{\Psi}_0 \mathbf{\tilde K}) \mathbf{a} \rvert ^2,
\end{equation}
where $\mathbf{\Psi}_0$ and $\mathbf{\Psi}_1 \in \mathbb{R}^{m \times M}$ can be written in the following form:
\begin{equation}
\mathbf{\Psi}_{0}=\left[\begin{array}{ccc}
\psi_{1}\left(\mathbf{u}_{0}\right) & \cdots & \psi_{M}\left(\mathbf{u}_{0}\right) \\
\psi_{1}\left(\mathbf{u}_{1}\right) & \cdots & \psi_{M}\left(\mathbf{u}_{1}\right) \\
\vdots & & \vdots \\
\psi_{1}\left(\mathbf{u}_{m-1}\right) & \cdots & \psi_{M}\left(\mathbf{u}_{m-1}\right)
\end{array}\right], \quad
\mathbf{\Psi}_{1}=\left[\begin{array}{ccc}
\psi_{1}\left(\mathbf{u}_{1}\right) & \cdots & \psi_{M}\left(\mathbf{u}_{1}\right) \\
\psi_{1}\left(\mathbf{u}_{2}\right) & \cdots & \psi_{M}\left(\mathbf{u}_{2}\right) \\
\vdots & & \vdots \\
\psi_{1}\left(\mathbf{u}_{m}\right) & \cdots & \psi_{M}\left(\mathbf{u}_{m}\right)
\end{array}\right].
\end{equation}

After the optimization, the operator $\mathbf{\tilde K}$ can be determined by
\begin{equation}
\mathbf{\tilde K}=\mathbf{\Psi}_0^{\dagger} \mathbf{\Psi}_1,
\end{equation}

The eigenfunctions of the Koopman operator and the Koopman modes then can be computed from the right eigenvectors and left eigenvectors of $\mathbf{\tilde K} $, respectively. The eigenvalues of $\mathbf{\tilde K}$ are approximations of the eigenvalues of $\mathbf{\mathcal{K}}$. For detailed derivation see \cite{morWilKR15}. A computational issue with EDMD is the expensive cost of computing the eigendecomposition of $\mathbf{\tilde K} \in  \mathbb{R}^{M \times M}$ when $M \gg m$, which is often the case in many applications.

Kernel DMD is proposed in \cite{morWilRK15} to reduce the computational cost of EDMD. This is done by using the compact SVD of the matrix $\mathbf{\Psi}_0= \mathbf{Q \mathbf{\Sigma} Z^T},  \mathbf{Q},  \mathbf{\Sigma} \in \mathbb R^{m\times m},  \mathbf{Z} \in \mathbb R^{M\times m}$. It is then proved in \cite{morWilRK15} that $\mathbf{\tilde K}$ has the same eigenvalues as the smaller matrix  $\mathbf{\hat K}=(\mathbf{\Sigma}^{-1}\mathbf{Q}^T) (\mathbf{\Psi}_1\mathbf{\Psi}_0^T) (\mathbf{Q} \mathbf{\Sigma}^{-1}) \in \mathbb{R}^{m \times m}$. Any right eigenvector $v$ of $\mathbf{\tilde K}$ corresponding to an eigenvalue $\lambda$ can be computed from the right eigenvector $\mathbf{\hat v} $ of $ \mathbf{\hat K}$ by $\mathbf{v} = \mathbf{Z} \mathbf{\hat v}$. From the SVD of  $\mathbf{\Psi}_0$, it is noticed that the eigendecomposition of $\mathbf{\Psi}_0 \mathbf{\Psi}_0^T$ is, 
\begin{equation}
\label{eq:eighatG}
\mathbf{\Psi}_0 \mathbf{\Psi}_0^T= \mathbf{Q} \mathbf{\Sigma}^2 \mathbf{Q}^T \in \mathbb R^{m\times m}.
\end{equation}
Therefore, if we can compute the eigendecomposition of $\mathbf{\Psi}_0 \mathbf{\Psi}_0^T$ and get  $\mathbf{Q}$, $\mathbf{\Sigma}$, then $\mathbf{\hat K}$ can be derived without SVD of $\mathbf{\Psi}_0$. The eigendecomposition of $\mathbf{\Psi}_0 \mathbf{\Psi}_0^T$ is of complexity $\mathcal{O}(m^3)$, which is less than $\mathcal{O}(M m^2)$, the SVD cost of $\mathbf{\Psi}_0 \mathbf{\Psi}_0^T$. It is further noticed that computing $\mathbf{\Psi}_0 \mathbf{\Psi}_0^T$ and $\mathbf{\Psi}_1 \mathbf{\Psi}_0^T$ is essentially implementing inner products of the two vectors $\psi(\mathbf{u}_i):=[\psi_1(\mathbf{u}_i) \quad \ldots \quad \psi_M(\mathbf{u}_i)]$ and $\psi(\mathbf{u}_j):=[\psi_1(\mathbf{u}_j) \quad \ldots \quad \psi_M(\mathbf{u}_j)], i, j=0, 1, \ldots, m$.  When $M$ is large, the computational cost of these inner products cannot be neglected. Usually, the observables include both the state variables and functions of them, making even $M \gg n$. The kernel function is then used to compute these inner products. As a result, the inner products in $\mathbb R^M$ are equivalently transformed to inner products in $\mathbb R^n$. This reduces the computations of directly computing the inner products $\psi(\mathbf{u}_i) \psi(\mathbf{u}_j)^T$. Please refer to \cite{morWilRK15} for a detailed explanation using illustrative examples. The final kernel DMD algorithm is reviewed in \Cref{algorithm:kerneldmd}, where $f(\mathbf{u}_i, \mathbf{u}_j)$ is a kernel function. Since $\mathbf{Z}$ can be represented by $\mathbf{\Psi}_0$, $\mathbf{Q}$ and $\mathbf{\Sigma}$, the eigenvectors of $\mathbf{\tilde K}$ can also be recovered without SVD of $\mathbf{\Psi}_0$. Furthermore, the eigenmodes of $\mathbf{\tilde K}$ are also computed independently of the SVD of $\mathbf{\Psi}_0$, see Step 7 in \Cref{algorithm:kerneldmd}. For a detailed derivation of it, please refer to \cite{morWilRK15}. Some common kernel functions that can be chosen are the polynomial kernel $f(\mathbf{x},\mathbf{y}) = (1+\mathbf{y}^T \mathbf{x})^{\alpha}$ or Gaussian kernel $f(\mathbf{x},\mathbf{y}) = \exp \left(-\left\|\mathbf{x}-\mathbf{y}\right\|^{2}/\sigma^2\right)$.

\begin{algorithm}[h]
	\begin{algorithmic}[1]
        \caption{Kernel DMD \cite{morWilRK15}} 
        \label{algorithm:kerneldmd}
        \State Compute the elements of $ \mathbf{\Psi}_0 \mathbf{\Psi}_0^T$ and $ \mathbf{\Psi}_1 \mathbf{\Psi}_0^T$ by kernel function: $ (\mathbf{\Psi}_0 \mathbf{\Psi}_0^T)_{ij} = f(\mathbf{u}_i, \mathbf{u}_j)$ and $ (\mathbf{\Psi}_1 \mathbf{\Psi}_0^T) =  f(\mathbf{u}_{i+1}, \mathbf{u}_j)$, with $i,j = 0 , \ldots, m-1$.
        \State Compute the eigendecomposition of $\mathbf{\Psi}_0 \mathbf{\Psi}_0^T$ via \cref{eq:eighatG} to get $\mathbf{Q}$ and $\mathbf{\Sigma}$.
        \State (optional) Choose the truncation rank $r$ that is smaller than the rank of $\mathbf{\Psi}_0 \mathbf{\Psi}_0^T$ to achieve a further reduction of the computation. Truncate the matrices $\mathbf{Q}$, and $\mathbf{\Sigma}$ by keeping the first $r$ columns of $\mathbf{Q}$ and first $r$ diagonal blocks $\mathbf{\Sigma}$ to obtain $\mathbf{Q}_r$ and $\mathbf{\Sigma}_r$.
        \State Compute $\mathbf{\hat K} = (\mathbf{\Sigma}_{r}^{-1}\mathbf{Q}_{r}^T) (\mathbf{\Psi}_1 \mathbf{\Psi}_0^T) (\mathbf{Q}_{r} \mathbf{\Sigma}_{r}^{-1})$.
        \State Compute the eigendecomposition of $\mathbf{\hat K} \mathbf{\hat W} = \mathbf{\hat W} \mathbf{\hat \Lambda}$ with $\mathbf{\Lambda} = diag(\lambda_1,\ldots,\lambda_r)$.
        \State Compute the approximated Koopman eigenfuction $\phi_k(\mathbf{u}) = (\mathbf{\psi}(\mathbf{u}) \mathbf{\Psi}_0^T)(\mathbf{Q}_{r} \mathbf{\Sigma}_{r}^{-1} \mathbf{\hat w}_k)$, where the inner product $\mathbf{\psi}(\mathbf{u}) \mathbf{\Psi}_0^T$ is obtained by evaluating the kernel function as explained in Step 1. $\mathbf{\hat w}_k$ is the right eigenvector from $\mathbf{\hat W} = [\mathbf{\hat w}_1, \cdots, \mathbf{\hat w}_r]$. 
        \State Set the Koopman modes as $\mathbf{v}_k =  (\hat \xi_k^* \mathbf{\Sigma}_r^{-1} \mathbf{Q}^{T}_{r} \mathbf{X}_0)^T \in \mathbb R^n$, where $\mathbf{\hat \xi}_k$ is the left eigenvector of the matrix $\mathbf{\hat K}$, and $\mathbf{\hat \xi}_k^* \mathbf{\hat w}_k = 1$ with $k=1, \ldots, r$.
        \State With eigenvalues, eigenfunctions and Koopman modes $\lambda_k, \phi_k, \mathbf{v}_k$, the approximation of the evolution can be done via \cref{eq:exdmd_F2K}.
	\end{algorithmic}
\end{algorithm}

However, either exact DMD or extended/kernel DMD cannot be straightforwardly applied to parametric problems, where the solution depends not only on the initial solution but also on the parameter variations. The parametric behavior of the solution usually cannot be captured by the DMD modes corresponding to any fixed value of the parameter provided by the DMD method. In the next section, we extend DMD to parametric DMD based on the RBF network.

\section{Proposed Parametric DMD}
\label{sec:work}

In many applications, parametric systems are widely used in multi-query tasks, such as optimal design, control, or uncertainty quantification. In this work, we consider parametric systems in a general form as,
\begin{equation}
\label{eq:eq_pDMD}
\begin{array}{ll}
\frac{d\mathbf{u}(t,\mu)}{dt}=g(\mathbf{u}(t,\mu), \mu), \quad
\mathbf{u}(t_0)=\mathbf{u}_0,\\
\mathbf{y}(t, \mu) = s(\mathbf{u}(t, \mu)).
\end{array}
\end{equation}
where $\mu \in \mathbb{R}^{n_p}$ is the vector of parameters, $\mathbf{u}(t, \mu) \in \mathbb{R}^{n}$ is the vector of states, and $\mathbf{y}(t,\mu) \in \mathbb{R}^{n_0}$ is the quantity of interest, also called the output. Existing DMD methods can not compute DMD modes which are also parametric, and as a result, they can only reconstruct the dynamics corresponding to a fixed value of $\mu$. Whenever the parameter value changes, DMD has to be reimplemented from scratch. In this section, we propose combining DMD with the RBF network to construct non-intrusive ROMs for parametric systems, which can predict the system's dynamics in both the parameter domain and the time domain. In \Cref{subsec:RBF}, we first review the RBF network, then in \Cref{subsec:DMD-RBF}, we connect it with DMD to realize parametric DMD.

\subsection{Radial Basis Function Network}%
\label{subsec:RBF}

The RBF method uses the weighted kernel function to approximate a given function $f(x)$: $ \mathbb{R}^{n_p} \to \mathbb{R}$ based on the data of $f(x)$. The approximate function $\hat{f}(x)$ constructed by RBF can be written as weighted summation of the RBFs, i.e.,
 
\begin{equation}\label{eq:RBF}
f(x) \approx \hat f(x) = \sum_{i=1}^{n_p} w_{i} \kappa \left(\left\|x-x_{i}\right\|\right).
\end{equation}

The kernel function is radially symmetric based on Euclidean distance $\left\|x-x_{i}\right\|$ or comparable metrics. The coefficients or weights $w_{i}, i=1,\ldots, n_p$ are determined by solving the following linear system of equations:
\begin{equation}\label{eq:rbf_weights}
\left[\begin{array}{ccc}
\kappa \left( \big\lVert x_1-x_1 \big\rVert \right) & \ldots & \kappa \left(\left\|x_{1}-x_{n_p}\right\|\right) \\
\vdots & \vdots & \vdots \\
\kappa \left(\left\|x_{n_p}-x_{1}\right\|\right) & \ldots & \kappa \left(\left\|x_{n_p}-x_{n_p}\right\|\right)
\end{array}\right]\left[\begin{array}{c}
w_{1} \\
\vdots \\
w_{n_p}
\end{array}\right]=\left[\begin{array}{c}
f\left(x_{1}\right) \\
\vdots \\
f\left(x_{n_p}\right)
\end{array}\right].
\end{equation}

\begin{figure}
\centering{
\begin{circuitikz}[scale=1.1, transform shape, american voltages, american currents]

\draw (0, 0) circle [radius = 0.45] node {$x_{n_p}$};

\draw (0, 1.1) circle [radius = 0.45] node {$x_{n_p-1}$};

\draw (0, 2.05) node {$\vdots$};

\draw (0, 2.8) circle [radius = 0.45] node {$x_2$};

\draw (0, 3.9) circle [radius = 0.45] node {$x_1$};

\draw (4, 0) circle [radius = 0.45] node {$\kappa_{n_p}$};

\draw (4, 1.1) circle [radius = 0.45] node {$\kappa_{n_p-1}$};

\draw (4, 2.05) node {$\vdots$};

\draw (4, 2.8) circle [radius = 0.45] node {$\kappa_2$};

\draw (4, 3.9) circle [radius = 0.45] node {$\kappa_1$};

\draw[-latex] (0.45, 0) to (3.55,0) {};

\draw[-latex] (0.45, 1.1) to (3.55,1.1) {};

\draw[-latex] (0.45, 2.8) to (3.55,2.8) {};

\draw[-latex] (0.45, 3.9) to (3.55,3.9) {};

\draw[-latex] (0.45, 0) to (3.55,1) {};

\draw[-latex] (0.45, 0) to (3.55,2.6) {};

\draw[-latex] (0.45, 0) to (3.55,3.6) {};

\draw[-latex] (0.45, 1.1) to (3.55,0.1) {};

\draw[-latex] (0.45, 1.1) to (3.55,2.7) {};

\draw[-latex] (0.45, 1.1) to (3.55,3.7) {};

\draw[-latex] (0.45, 2.8) to (3.55,3.8) {};

\draw[-latex] (0.45, 2.8) to (3.55,1.2) {};

\draw[-latex] (0.45, 2.8) to (3.55,0.2) {};

\draw[-latex] (0.45, 3.9) to (3.55,2.9) {};

\draw[-latex] (0.45, 3.9) to (3.55,1.3) {};

\draw[-latex] (0.45, 3.9) to (3.55,0.3) {};

\draw[-latex] (4.45, 0) to (7.55,1.75) {};

\draw[-latex] (4.45, 1.1) to (7.55,1.85) {};

\draw[-latex] (4.45, 2.8) to (7.55,1.95) {};

\draw[-latex] (4.45, 3.9) to (7.55,2.05) {};

\draw (8, 1.95) circle [radius = 0.45] node {$\hat{f}(x)$};

\draw (5.3, 0.76) node {$w_{n_p}$};

\draw (5.3, 1.55) node {$w_{n_p-1}$};

\draw (5.3, 2.75) node {$w_2$};

\draw (5.3, 3.62) node {$w_1$};

\end{circuitikz}}
\caption{RBF interpolation as a shallow neural network.}\label{fig:RBF_snn}
\end{figure}
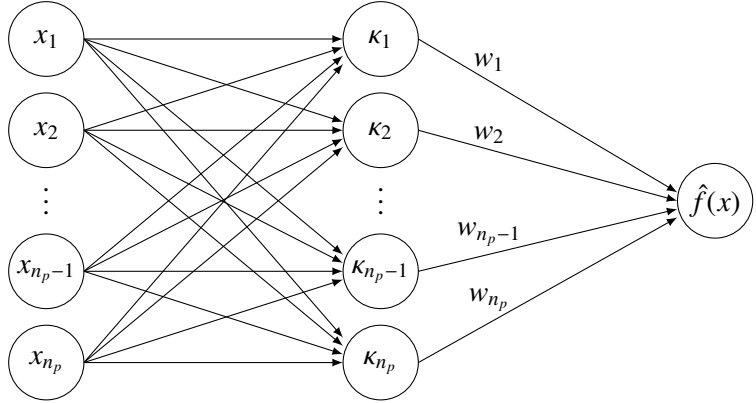

The whole process of computing $\hat f(x)$ works like a shallow neural network shown in \Cref{fig:RBF_snn}, where $\kappa_i = \kappa \left(\left\|x-x_{i}\right\|\right)$. After the weights are fixed, the interpolation can be completed simply using the weighted summation in \cref{eq:RBF}. The detailed process of computing the approximate function $\hat f(x)$ is presented in \Cref{algorithm:rbf}.

In this work, when training the RBF network, the data points $x_i$ are the samples of the parameters $\mu_i$, and $f(x)$ corresponds to each entry of the solution vector $u(t, \mu)$ at any time instance $t_k \leq T_0$ and any training sample of $\mu$, i.e., each entry in the snapshot matrices $\mathbf{X}_0$ and $\mathbf{X}_1$ in \cref{eq:x0x1}.

\begin{algorithm}[h]
    \begin{algorithmic}[1]
        \caption{RBF network construction} \label{algorithm:rbf}
        \State Choose the appropriate kernel function $\kappa$ and its shape factor $\varepsilon$ when necessary.
        \State Compute the Euclidean distance $r = \left\|x-x_{i}\right\|$ between data points.
        \State Compute the radial basis function $\kappa(r)$ and construct the coefficient matrix in \cref{eq:rbf_weights}.
        \State Determine the coefficients $w_i$ by solving the linear system in \cref{eq:rbf_weights}.
    \end{algorithmic}
\end{algorithm}

The kernel functions can be chosen in a wide variety, such as splines, Gaussian, Multi-quadrics, and so on. \Cref{tab:kernels} provides a chart with some commonly used basis functions. In this work, inverse multi-quadrics (IMQ) is used with shape factor $\varepsilon = 1/30$.

\begin{table}[h]
\begin{center}
\begin{tabular}{|l|l|}
\hline Linear splines & $\left\|x-x_{i}\right\|$ \\
\hline Cubic splines & $\left\|x-x_{i}\right\|^{3}$ \\
\hline Thin plate splines & $\left\|x-x_{i}\right\|^{k} \ln \left\|x-x_{i}\right\| ; k \in[2,4, \ldots]$ \\
\hline Multi-quadrics & $\left(1+\varepsilon^{2}\left\|x-x_{i}\right\|^{2}\right)^{1 / 2}$ \\
\hline Inverse multiquadrics & $\left(1+\varepsilon^{2}\left\|x-x_{i}\right\|^{2}\right)^{-1 / 2}$ \\
\hline Gaussian & $\exp \left(-\varepsilon^{2}\left\|x-x_{i}\right\|^{2}\right)$ \\
\hline
\end{tabular}
\caption{Some commonly used RBFs, $x_{\mathrm{i}}$ denotes the $i$-th center of the RBF.}\label{tab:kernels}
\end{center}
\end{table}

\subsection{Parametric DMD framework}%
\label{subsec:DMD-RBF}

\begin{figure}
	\centering
	\includegraphics[width=1.0\linewidth]{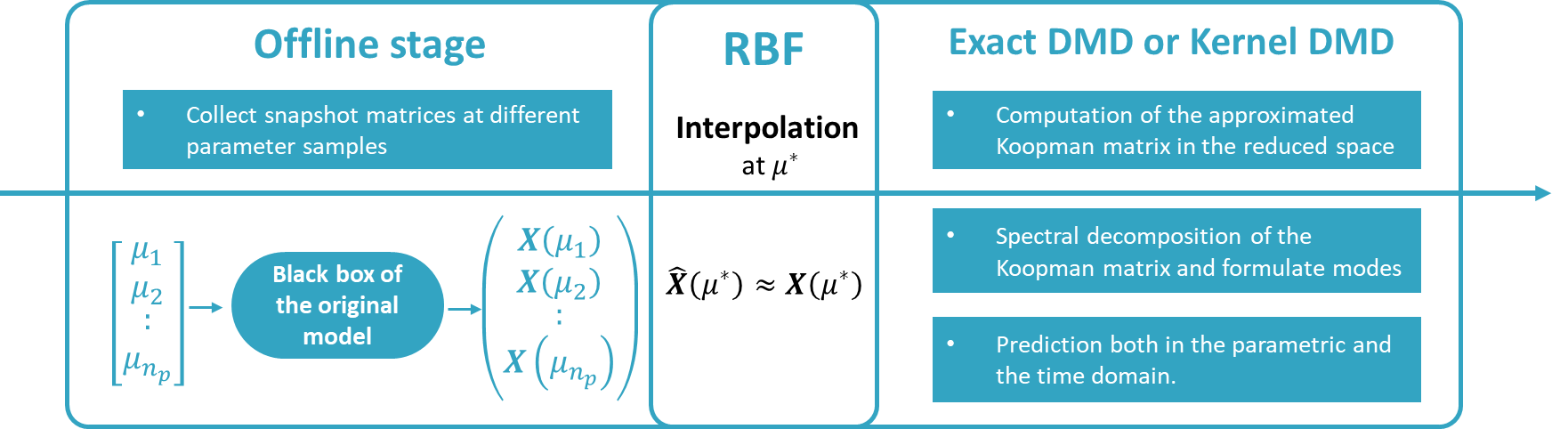}
	\caption{Flowchart of the parametric DMD algorithm.}
	\label{fig:workflow}
\end{figure}

In this section, we propose the parametric DMD framework. After the collection of snapshots at limited samples of training parameters,  the RBF network is first trained using these snapshots. The trained RBF network can then predict snapshots at any new parameter. After the snapshot matrices $\mathbf{X}_0, \mathbf{X}_1$ corresponding to the new parameter are computed, DMD is implemented on the new snapshot matrices to generate the DMD modes for predicting the solution in the time domain. The whole flow chart of the parametric DMD framework can be seen in \Cref{fig:workflow}. At the offline stage, the snapshot matrices $\mathbf{X}(\mu_k):=[\mathbf{u}_1(\mu_k), \ldots, \mathbf{u}_{m+1}(\mu_k)]$ corresponding to different samples $\mu_k, k=1,\ldots,n_p$, of the parameter $\mu$ are first computed via, e.g., black-box simulation of a dynamical system. These are used as training data for the RBF network. Then the RBF network is used to construct an approximate function $\mathbf{\hat X}_{ij}(\mu): \mu \mapsto \mathbb{R}$ for each entry  $\mathbf{X}_{ij}(\mu), i=1,\ldots, n, j=1,\ldots, m+1$, of a snapshot matrix function $\mathbf{X}(\mu)$. More specifically, $\hat f(x)$ in \cref{eq:RBF} now becomes $\mathbf{\hat X}_{ij}(\mu)$, and $f(x)$ is now $\mathbf{X}_{ij}(\mu)$. The RBF network is used to learn the $i,j$-th entry of $\mathbf{X}(\mu)$ using the data $\mathbf{X}_{ij} (\mu_k)$, i.e., the $i,j$-th entry of the snapshot matrices $\mathbf{X}(\mu_k)$ at the $n_p$ parameter samples $\mu_k, k=1,\ldots,n_p$. The weights $w_i$ in \cref{eq:RBF} are computed once for each entry  $\mathbf{X}_{ij}(\mu)$. After the weights are computed, the RBF network $ \mathbf{\hat X}_{ij}(\mu)$ for the $i,j$-th entry is trained and is ready to be used at the online stage. The predicted snapshot matrix at $\mu^*$ is nothing but $\mathbf{\hat X}(\mu^*)$.

At the online stage, instead of repeated black-box simulation of the large-scale model in \cref{eq:eq_pDMD}, the maps $\mathbf{\hat X}_{ij}(\mu)$ constructed by the RBF networks are called to compute the approximated snapshot matrix $\mathbf{\hat X}(\mu^*)$ at any new parameter sample $\mu^*$. $\mathbf{\hat X}(\mu^*)$ is then split into two snapshot matrices $\mathbf{\hat X}_0 \in \mathbb{R}^{n\times m}$ and $\mathbf{\hat X}_1 \in \mathbb R^{n\times m}$. For example, if $\mathbf{\hat X}(\mu^*) \in \mathbb R^{n\times m+1}$ approximates $\mathbf{X}(\mu^*):=[\mathbf{u}_1(\mu^*), \ldots, \mathbf{u}_{m+1}(\mu^*)]$, then $\mathbf{\hat X}_0(\mu^*)=\mathbf{\hat X}(\mu^*)[:,1\!:\!m]$, and $\mathbf{\hat X}_1(\mu^*)=\mathbf{\hat X}(\mu^*)[:, 2\!:\!m\!+\!1]$. Here we use the MATLAB notation for matrix blocks. The exact DMD or the kernel DMD is then applied to $\mathbf{\hat X}_0(\mu^*)$ and $\mathbf{\hat X}_1(\mu^*)$ to predict the time evolution of the solution corresponding to $\mu^*$. In summary, the RBF networks are used to predict the dynamics in the parameter domain and the DMD is employed for the time-evolution prediction. This process of parametric DMD is detailed in \Cref{algorithm:dmd-rbf}. 

\begin{algorithm}[h]
    \begin{algorithmic}[1]
        \caption{Parametric DMD framework} \label{algorithm:dmd-rbf}
        \State Sample (equidistant or random) the parameters in the range of the parameter space.
        \State Generate the snapshot matrices $\mathbf{X}(\mu_k)$ for each parameter sample $\mu_k$ with $k=1,\ldots,n_p$;
        \State Train the RBF network for each entry of $\mathbf{X}(\mu)$ by \Cref{algorithm:rbf}, and get the approximate snapshot matrix function $\mathbf{\hat X}(\mu)$ in the form of RBF networks.
        \State For any new parameter sample $\mu^*$, evaluate the approximate snapshot matrix function at $\mu^*$ and take $\mathbf{\hat X}(\mu^*)$ as the new snapshot matrix corresponding to $\mu^*$. Split it into $\mathbf{\hat X}_0(\mu^*)$ and $\mathbf{\hat X}_1(\mu^*)$.
        \State Apply the exact DMD (\Cref{algorithm:exactdmd}) or the kernel DMD (\Cref{algorithm:kerneldmd}) to $\mathbf{\hat X}_0(\mu^*)$ and $\mathbf{\hat X}_1(\mu^*)$ to reconstruct and predict the time-evolution of the dynamics corresponding to $\mu^*$.
    \end{algorithmic}
\end{algorithm}

\section{Numerical examples}
\label{sec:examples}

In this section, we test the performance of the proposed parametric DMD method with some models from engineering applications. Two examples are related to electrochemical processes. The first one considers lithium-ion battery model. Lithium-ion batteries are of high importance in the context of electromobility. Understanding of their dynamics is of high interest. The second example is a ferrocyanide reduction oxidation reaction. This is a common model system in electrochemistry which exemplifies diffusion controlled fast electrochemical process. As the last example, the FitzHugh-Nagumo model is a prototype of an excitable system, for example, a neuron. A common feature of all the models is that they are systems with parameters and time-varying inputs that can be considered as time-varying parameters. In the following subsections, we discuss the numerical tests on each of them separately. In all the figures illustrating the numerical results, ``RBF-DMD'' represents parametric DMD, and ``reference'' refers to the solution computed by directly simulating the original model. According to the error computation in the numerical examples, we use relative error at any testing parameter $\mu^*$ defined as follows:
\begin{equation}
\label{eq:rel_err}
\epsilon_i(t, \mu^*) = \frac{\left|y_i(t, \mu^*) - \hat y_i(t, \mu^*) \right|}{\max\limits_{0\leq t_j \leq T} \left|y_i(t_j, \mu^*)\right|}.
\end{equation}
Here the index $i$ means the $i$-th output, i.e., the $i$-th entry of $y(t, \mu^*) \in \mathbb{R}^{n_0}$. To evaluate the performance of the proposed method in the parameter domain, the time-average relative error is used and is defined as:
\begin{equation}
\label{eq:ta_rel_err}
\epsilon_i^{ave}(\mu^*) = \frac{1}{n_T} \sum_{j=0}^{n_T-1} \epsilon_i(t_j, \mu^*).
\end{equation}
As for the computation time, on the one hand, the snapshot generation and the RBF training are run only once at the offline stage. On the other hand, the runtimes of the RBF prediction, the DMD prediction and the FOM simulation at the online phase are respectively the average values of the runtimes over all the testing parameters.

\subsection{Lithium-ion Battery Model}%
\label{subsubsec:p2dpsd_battery}

\begin{figure}[ht]
\begin{subfigure}{0.5\textwidth}
\includegraphics[width=0.95\linewidth, height=3.5cm]{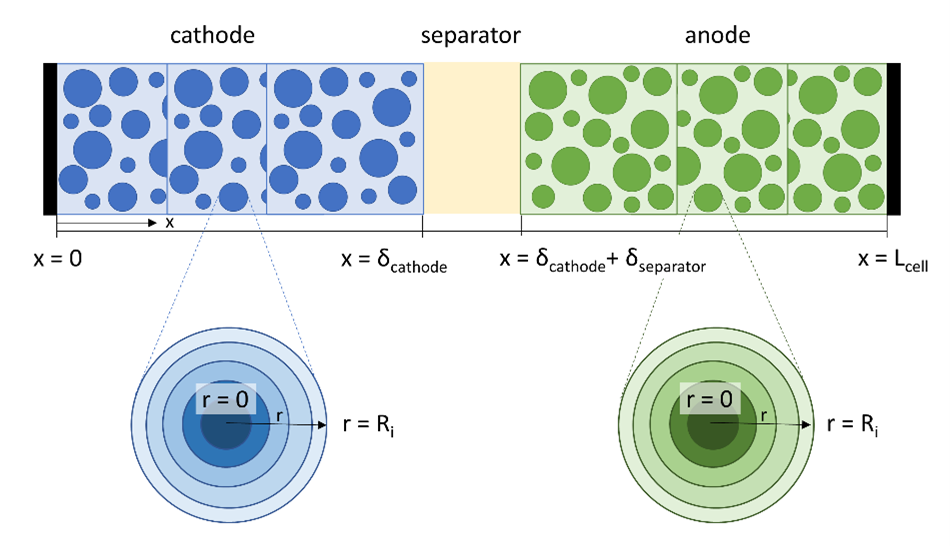} 
\caption{}
\label{fig:Battery_P2D}
\end{subfigure}
\begin{subfigure}{0.5\textwidth}
\centering
\includegraphics[width=0.7\linewidth, height=3.5cm]{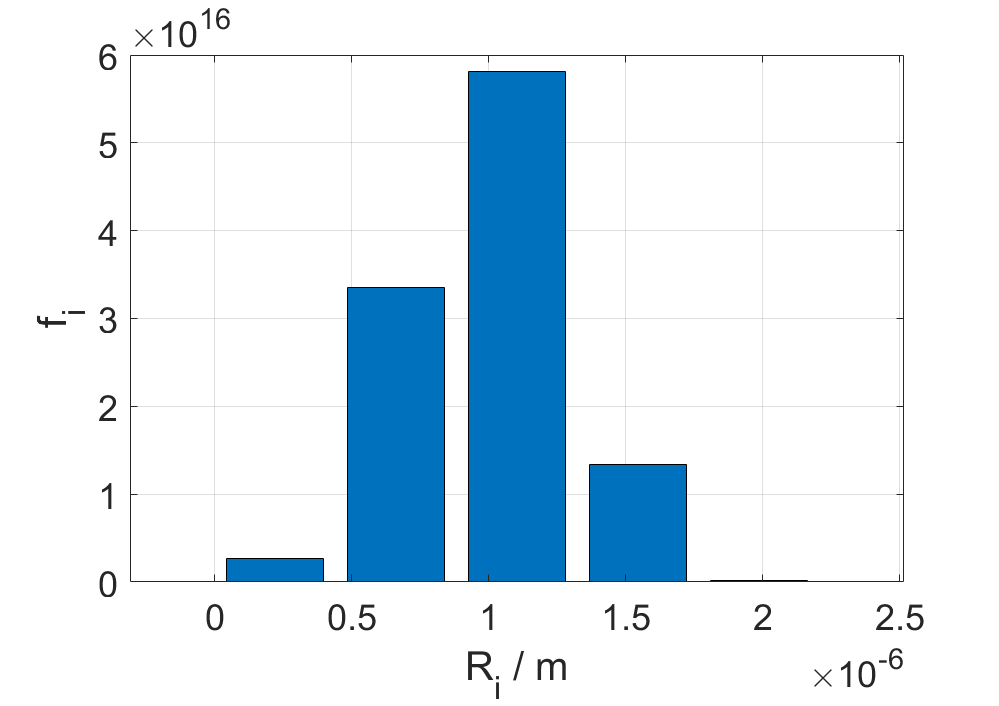}
\caption{}
\label{fig:Weibull}
\end{subfigure}
\caption{Mathematical model of the lithium-ion battery. (a) Schematic of the modelling domains of the battery model. (b) Particle size distribution of electrode (anode) assumed in P2D battery model simulation.}
\label{fig:Li-ion-battery}
\end{figure}

As an example for validating of the proposed methodology, we consider the widely implemented yet complicated mathematical model of a lithium-ion battery, the pseudo-two-dimensional (P2D) battery model, which was previously introduced in \cite{DoyFN93}. \Cref{fig:Battery_P2D} depicts the schematic of the P2D battery model. As the name suggests, the P2D battery model comprises two modelling scales: the computation of lithium concentration and potential gradients across the battery model (macro-scale) as well as the diffusion of lithium ions within the electrode (micro-scale).

Further complexities arise in the P2D battery model when one considers a distribution of different particle sizes in the electrode (anode), which has been introduced by Röder et al. \cite{RodSSetal16}. Considering different sizes of the solid particles, the battery dynamics at a wider operational condition can be better reproduced via a model-based approach. Here, we assume the particle size distribution within the electrode follows a Weibull distribution density defined as:
\begin{equation}
    h\left(R_{i} ; c_1, c_2\right)=c_1 c_2\left(c_1 R_{i}\right)^{c_2-1} e^{-\left(c_1 R_{i}\right)^{c_2}},
\end{equation}
where $c_1$ is the scaling factor and $c_2$ is the form factor of the distribution density. $R_i$ is the particle radius of the $i$-th particle size class in the electrode. \Cref{fig:Weibull} shows the simulated particle size distribution of the electrode with five different radius classes with $c_1=9.064 \times 10^5$ and $c_2=4$. Summing up the surface and volume densities across every particle radius class yields the surface area ratio as well as the volume fractions of the total active materials in the battery:

\begin{equation}
\begin{aligned}
a_{s} &=\int_{0}^{\infty} f_{\text {area }}\left(R_{i}\right) d R_{i}, \\
\varepsilon_{s} &=\int_{0}^{\infty} f_{\text {vol }}\left(R_{i}\right) d R_{i}.
\end{aligned}
\end{equation}

The governing equations for the P2D-PSD model are derived from the conservation laws of species and charge transport. The governing equations of the P2D battery model incorporating the effect of particle size distribution are detailed in \Cref{tab:P2D_eq}.

\begin{table}[htbp]
{\footnotesize
  \caption{Governing equations of the P2D battery model with the corresponding boundary conditions \cite{LegRKetal14}. Subscript $i$ describes the $i^{th}$ particle size classes, $s$ for solid phase, $e$ for electrolyte phase.}  \label{tab:P2D_eq}
\begin{center}
\resizebox{\textwidth}{!}{%
  \begin{tabular}{|c|c|} \hline
   Model equations & Boundary conditions \\ \hline
    $\frac{\partial c_{s}\left(r, R_{i}\right)}{\partial t}=\frac{1}{r^{2}} \frac{\partial}{\partial r}\left(D_{s, a/c} \cdot r^{2} \cdot \frac{\partial c_{s}\left(r, R_{i}\right)}{\partial r}\right)$ & \vtop{\hbox{\strut $\frac{\partial c_{s}\left(r = 0, R_{i}\right)}{\partial r} = 0$}\hbox{\strut $\frac{\partial c_{s}\left(r=R_{i}, R_{i}\right)}{\partial r}=\frac{j^{L i}\left(x, R_{i}\right)}{a_{s} D_{s} F}$}} \\ \hline
    $\varepsilon_{e} \frac{\partial C_{e}(x)}{\partial t}=\frac{\partial}{\partial x}\left(D_{e, eff} \cdot \frac{\partial C_{e}(x)}{\partial x}\right)+\left(1-t_{p}\right) \cdot \frac{j^{L_{i}}(x)}{F}$ & \vtop{\hbox{\strut $\frac{\partial c_{e}\left(x = 0\right)}{\partial x} = 0$}\hbox{\strut $\frac{\partial c_{e}\left(x = L_{cell}\right)}{\partial x}=0$}} \\ \hline
    \vtop{\hbox{\strut $J_{s}(x)=-\sigma_{s} \varepsilon_{s} \frac{\partial \phi_{s}(x)}{\partial x}$}
\hbox{\strut$\frac{\partial J_{s}(x)}{\partial x}=-\left(j^{L i}(x)+j^{D L}(x)\right)$ }}& \vtop{\hbox{\strut$\frac{\partial J_{s}(x=0)}{\partial x}=a_{s} \frac{I}{A_{\text {cell }}}$} \hbox{\strut $\frac{\partial J_{s}\left(x=L_{\text {cell }}\right)}{\partial x}=-a_{s} \frac{I}{A_{\text {cell }}}$}}\\ \hline
    \vtop{\hbox{\strut $J_{e}(x)=-\sigma_{e}(x) \frac{\varepsilon_{e}}{\tau} \frac{\partial \phi_{e}(x)}{\partial x}-\sigma_{D e} \frac{\varepsilon_{e}}{\tau} \frac{\partial \ln \left(c_{e}(x)\right)}{\partial x}$}\hbox{\strut $\frac{\partial J_{e}(x)}{\partial x}=j^{L i}(x)+j^{D L}(x)$}}& \vtop{\hbox{\strut $\phi_{e}(x=0) = 0$} \hbox{\strut $\frac{\partial \phi_{e}(x=L_{cell})}{\partial x} = 0$}}\\ \hline
    \multicolumn{2}{|c|}{$j^{D L}(x)=a_{s} C_{D L} \frac{\partial(\Delta \phi(x))}{\partial t}$}\\ \hline
    \multicolumn{2}{|c|}{\vtop{\hbox{\strut$j^{L i}\left(x, R_{i}\right)=a_{s, j}\left(R_{i}\right) j_{0}\left(R_{i}\right)\left(\exp \left(\frac{\alpha \eta\left(x, R_{i}\right) F}{R T}\right)-\exp \left(\frac{(1-\alpha) \eta\left(x, R_{i}\right) F}{R T}\right)\right)$}
\hbox{\strut $\eta\left(x, R_{i}\right)=\phi_{s}(x)-\phi_{e}(x)-U\left(x, R_{i}\right)$}
\hbox{\strut $j^{L i}(x)=\int j^{L i}\left(x, R_{i}\right) d R_{i}$}}}\\ \hline
  \end{tabular}}
\end{center}
}
\end{table}

\begin{table}
{\footnotesize
  \caption{Parameter set for P2D battery model simulation.}  \label{tab:P2D_parameters}
\begin{center}
  \begin{tabular}{|c|l|l|l|} \hline
  Symbol & Parameter & Unit & Value \\ \hline
  $R$ & Gas constant & $J mol^{-1} K^{-1}$ & $8.314$ \\ \hline
  $F$ & Faraday constant & $C mol^{-1}$ & $96485$ \\ \hline
  $T$ & Temperature & $K$ & $298$ \\ \hline
  $R_c$ & Radius of cathode	& $\mu m$ & $1$ \\ \hline
  $\delta_{anode}$	 & Anode’s thickness	& $\mu m$ & $50$ \\ \hline
  $\delta_{separator}$ & Separator’s thickness	& $\mu m$ & $26.4$ \\ \hline
  $\delta_{cathode}$ & Cathode’s thickness & $\mu m$ & $25.4$ \\ \hline
  $D_{s,anode}$ & Diffusion coefficient anode	 & $m^2 s^{-1}$	& $2\times 10^{-16}$ \\ \hline
  $D_{s,cathode}$ & Diffusion coefficient cathode	 & $m^2 s^{-1}$	& $3.7\times 10^{-16}$ \\ \hline
  $\tau_{anode}$ & Tortuosity anode	&-&	$3.67$ \\ \hline
  $\tau_{separator}$ & Tortuosity separator &-& $1.4$ \\ \hline
  $\tau_{cathode}$ & Tortuosity cathode	&-&	$3.67$ \\ \hline
  $\varepsilon_{e,anode}$ & Volume fraction electrolyte anode &-& $0.33$\\ \hline
  $\varepsilon_{e,cathode}$ & Volume fraction electrolyte cathode &-& $0.33$\\ \hline
  $\varepsilon_{e,separator}$ & Volume fraction electrolyte separator &-& $0.5$\\ \hline
  $\varepsilon_{s,cathode}$ & Volume fraction cathode &-& $0.5$\\ \hline
  $\alpha$ & Charge transfer coefficient &-& $0.5$ \\ \hline
  $C_{DL,anode}$ & Double layer capacitance anode & $F m^{-2}$& $0.2$ \\ \hline
  $C_{DL,cathode}$ & Double layer capacitance cathode & $F m^{-2}$& $0.2$ \\ \hline
  $t_p$ & Transference number  &-& $0.37$ \\ \hline
  $\sigma_{s, anode}$ & Electrical conductivity anode	&$S m^{-1}$&	$100$ \\ \hline
  $\sigma_{s, cathode}$ & Electrical conductivity cathode	&$S m^{-1}$&	$10$ \\ \hline  
\end{tabular}
\end{center}}
\end{table}

\begin{figure}[ht]
	\centering
	\includegraphics[width=0.65\linewidth]{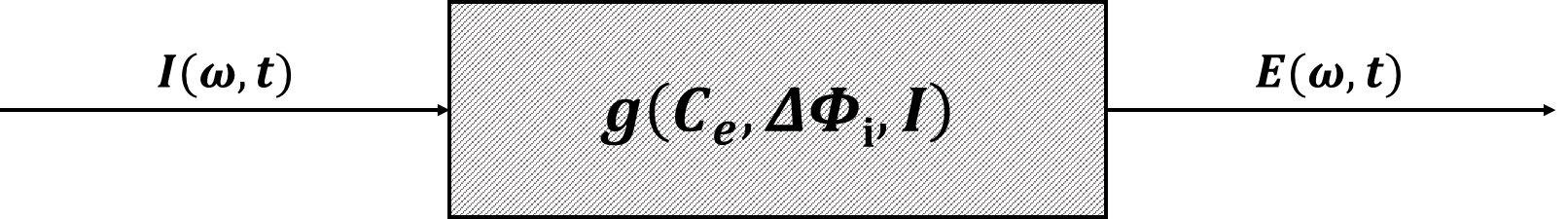}
	\caption{Graphical model of a Lithium-ion battery.}
	\label{fig:inout}
\end{figure}

It is also seen that all the governing equations are coupled with each other. Due to the high complexity of the coupled governing equations, it is almost impossible to extract the discretized system matrices and nonlinear terms from the spatial discretization of the PDEs given all the parameters are fixed, not to mention their parametrized forms. Consequently, projection-based MOR methods cannot be applied for MOR of this model, and the non-intrusive MOR is the only possible choice. That leads to the application of our proposed parametric DMD to this example. The input $I(\omega, t)$ is the current with a certain frequency $\omega \in [10^{-2},10^{4}]$ and the output $E(\omega, t)$ is the voltage, which can be shown as the difference of potential at the current collectors between the anode and the cathode. Both are shown in \cref{eq:inout}. The whole general in-output model is shown in \Cref{fig:inout}.
\begin{equation}
\label{eq:inout}
\begin{gathered}
I(\omega, t)=\frac{I}{A_{cell}} \sin (\omega t) \\
E(\omega, t)=\phi_{s}(x=0)-\phi_{s}\left(x=L_{cell}\right)
\end{gathered}
\end{equation}

The original spatially discretized ODE model has $n=325$ states. We use 100 snapshot matrices corresponding to 100 frequency training samples $\omega_i, i=1,\ldots, 100$, with 10-base logarithmic spacing in $[-2 ,4]$. At the offline stage, the RBF network is trained with these parameters. The snapshots corresponding to each frequency sample in a limited time interval $[0, T_0]$ are computed by an ODE solver: ode15s in MATLAB.  Here $T_0 =  T/2$, with $T$ being the final simulation time. That means the original model is simulated till half of the final simulation time to get the snapshot matrices. The dynamics corresponding to time span $[0, T]$ at any testing frequency will be predicted. At the online stage, the snapshot matrix function $\mathbf{\hat X}(\omega)$ in the form of RBF networks is evaluated at a new frequency sample $\omega^*$ to get an approximate snapshot matrix $\mathbf{\hat X}(\omega^*)$  that is considered as the new snapshot matrix. DMD is then applied to $\mathbf{\hat X}(\omega^*)$ to predict the output voltage $E(\omega, t)$ at any time $t>T_0$.

\begin{figure}
\begin{subfigure}{0.5\textwidth}
\includegraphics[width=0.95\linewidth]{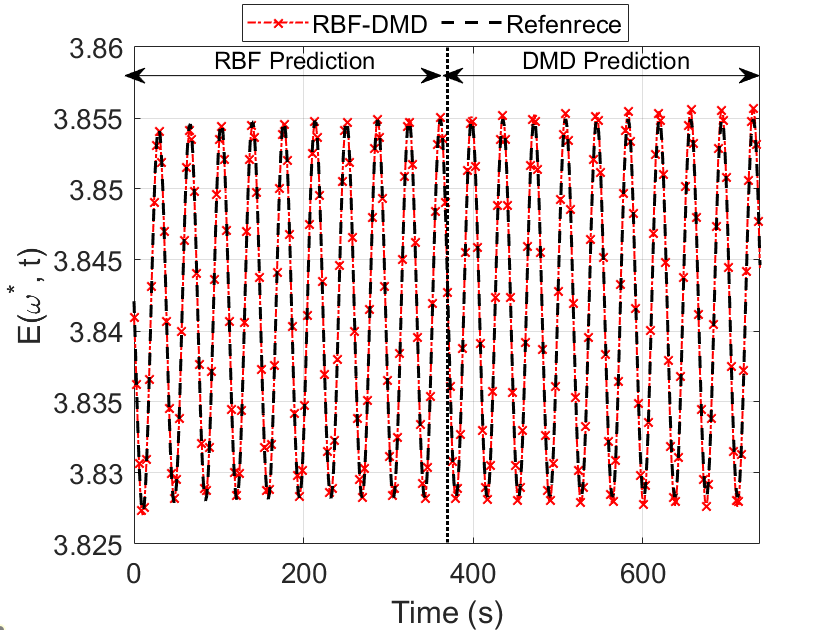} 
\caption{}
\label{fig:P2D_PSD_last_1}
\end{subfigure}
\begin{subfigure}{0.5\textwidth}
\includegraphics[width=0.95\linewidth]{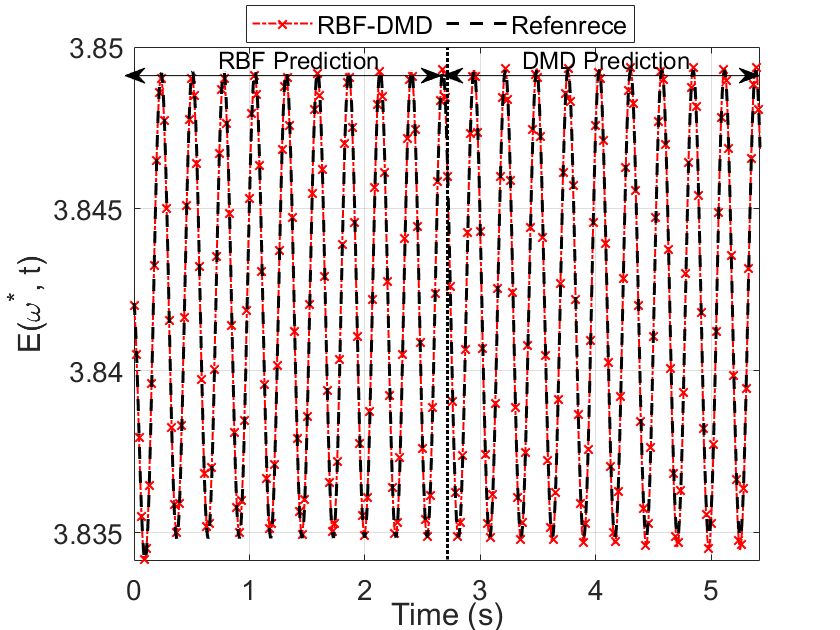}
\caption{}
\label{fig:P2D_PSD_last_2}
\end{subfigure}
\caption{Lithium-ion battery model: the parametric DMD solution vs the reference solution for $E(\omega^*, t)$. (a) $ \omega^* = 0.025 \,Hz$. (b) $ \omega^* = 3.69 \,Hz$.}
\label{fig:P2D_PSD_last}
\end{figure}

\begin{figure}
\begin{subfigure}{0.5\textwidth}
\includegraphics[width=0.95\linewidth]{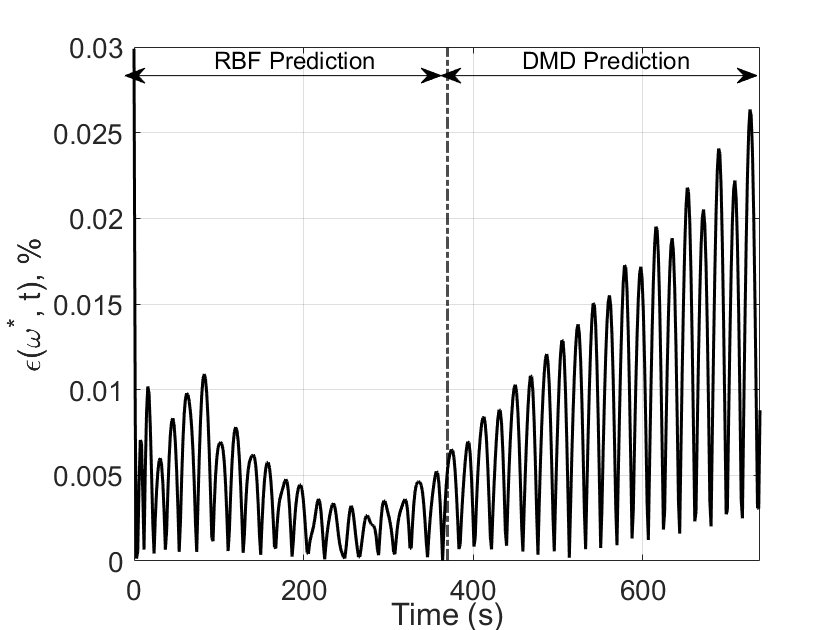} 
\caption{}
\label{fig:P2D_PSD_error_1}
\end{subfigure}
\begin{subfigure}{0.5\textwidth}
\includegraphics[width=0.95\linewidth]{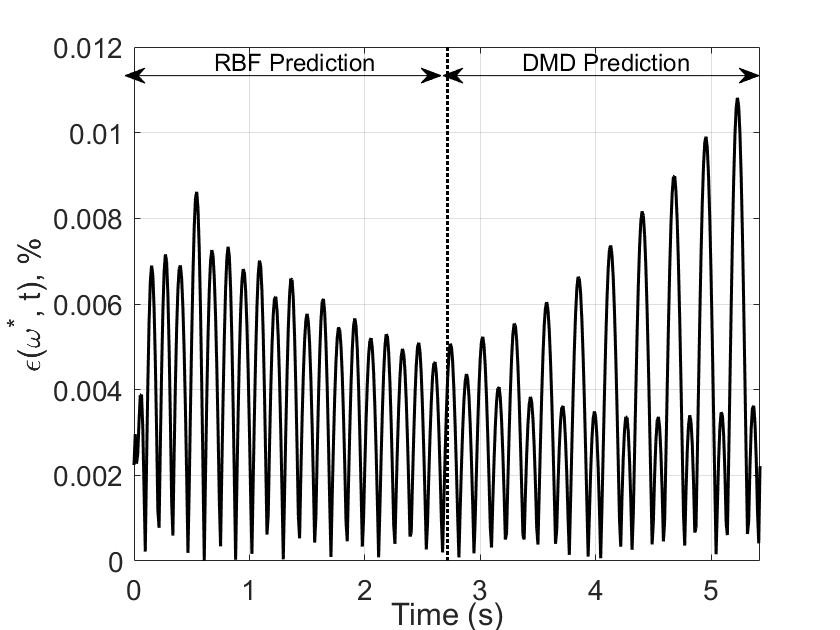}
\caption{}
\label{fig:P2D_PSD_error_2}
\end{subfigure}
\caption{Lithium-ion battery model: the relative error of the parametric DMD solution for $E(\omega^*, t)$. (a) $ \omega^* = 0.025 \,Hz$. (b)
$ \omega^* = 3.69 \,Hz$.}
\label{fig:P2D_PSD_error}
\end{figure}

\begin{figure}
\centering
\includegraphics[width=0.5\linewidth]{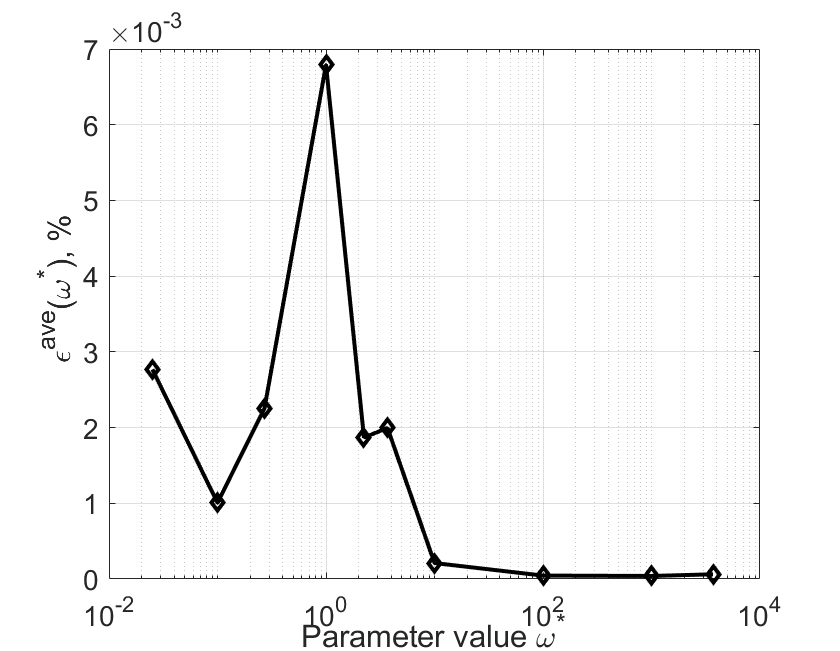}
\caption{Lithium-ion battery model: the time-average relative error of the parametric DMD for $E(\omega^*, t)$ at different testing frequencies $\omega^*$.}
\label{fig:P2D_PSD_merror}
\end{figure}

Exact DMD is employed in the proposed parametric DMD for this model. The results are derived by 14 dominant DMD modes, i.e., $r =14$ in \Cref{algorithm:exactdmd}. The time-evolution of the output voltage at $\omega^*$ computed by the ODE solver is considered as the reference solution. Both the reference solution and the output computed by the parametric DMD are presented in \Cref{fig:P2D_PSD_last}. The RBF-DMD solution is the voltage derived by the proposed parametric DMD. The voltage in $[0, T_0]$ is predicted by the RBF network. Based on this, DMD then predicts the evolution in $[T_0, T]$. The relative error between the reference voltage and the RBF-DMD voltage is presented in \Cref{fig:P2D_PSD_error}. It can be observed that the maximal relative error is under $0.03 \% $. The plot for the time-average relative error at different testing frequencies $\omega^* = 0.025, 0.1, 0.271, 1, 2.239, 3.690,10,100,1000, 3689.776 \,Hz$ is shown in \Cref{fig:P2D_PSD_merror}. RBF-DMD predicts the voltage with no more than $0.007 \% $ relative error compared to the reference solution both in low and high testing frequencies.
 
\begin{figure}
\centering
\includegraphics[width=0.9\linewidth]{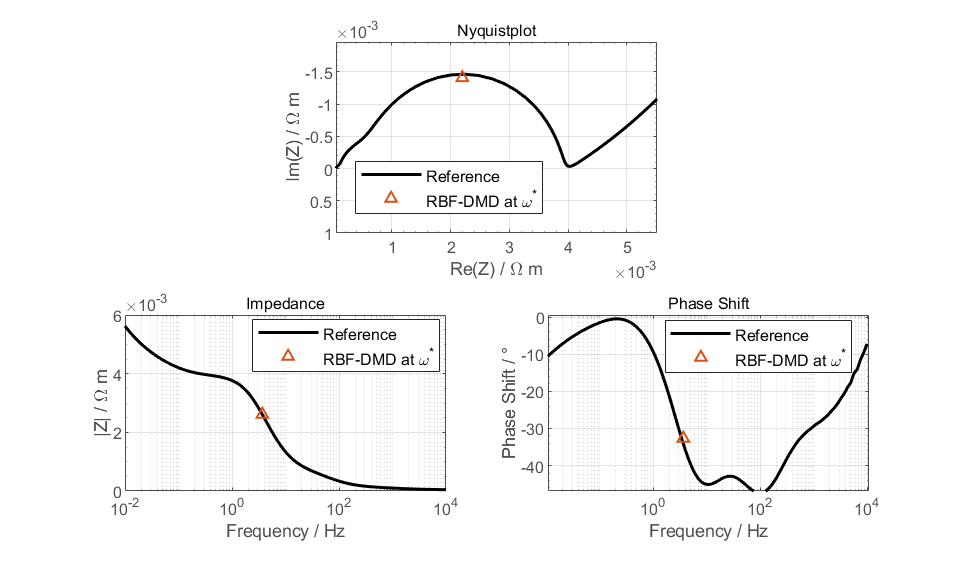}
\caption{Lithium-ion battery model: Nyquist and bode representations of the parametric DMD solution at $\omega^* =  3.69 \,Hz$ vs the reference solution.}
\label{fig:Bode}
\end{figure}

Electrochemical impedance spectroscopy (EIS) is commonly used to monitor the performance of the lithium-ion battery. When the input is the current with different frequencies, the output voltage is transformed from time to frequency domain by Fast Fourier Transformation (FFT) to analyse the model. In this example, the results of EIS are shown in Nyquist and Bode diagrams, see \Cref{fig:Bode}. The subfigure above is the Nyquist plot presenting the imaginary part of the complex impedance as a function of its real part. It can be observed that there exists a semicircle at the high frequency range and a non-vertical line at the intermediate frequency range, which can be interpreted as the resistance of the electrolyte and the resistance of the diffusive layer, respectively, in the practical application. In this subfigure, the complex impedance computed by parametric DMD at the new frequency $\omega^* = 3.69 \,Hz$ conforms to the pattern from the reference solution. The bottom-left one shows the relationship between the impedance and the frequency, while the bottom-right one is the phase shift changing with the frequency. Both subfigures show the great matching between the solution from the parametric DMD and the reference solution at testing frequency $\omega^*$. We can conclude that the proposed method delivers satisfactory accuracy in the parameter space and the time domain.


The runtime comparison for this example is shown in \Cref{tab:P2D_computation_time}. The computation time of parametric DMD includes the offline stage and the online stage. The offline stage of computing all the snapshots at 100 training samples takes 171.581 seconds. Training the RBF network at the offline stage takes 0.445 seconds. The online RBF prediction at a new parameter sample in the training time interval costs 0.007 seconds. The online DMD prediction in the future time interval takes 1.697 seconds. Computing the reference solution at one testing sample of $\omega^*$ via ODE solver i.e., the FOM simulation needs 4.787 seconds. The online speed-up is around $2-3$ times faster. It is clear that if the original model needs to be simulated to get the output response at more than 40 different values of $\omega$, the proposed parametric DMD method will outperform the direct simulation without MOR.

\begin{table}[h]
{\footnotesize
  \caption{Lithium-ion battery model: The computation time (seconds) of parametric DMD and that of the FOM simulation.}  \label{tab:P2D_computation_time}
\begin{center}
  \begin{tabular}{|c|c|c|c|c|} \hline
  Snapshot generation & RBF training & RBF prediction & DMD prediction & FOM simulation\\ \hline
  171.581 &0.445  & 0.007 & 1.697 & 4.787 \\ \hline
\end{tabular}
\end{center}}
\end{table}

\subsection{Coupled electrochemical kinetics and diffusion model}%
\begin{figure}
	\centering
	\includegraphics[width=0.8\linewidth]{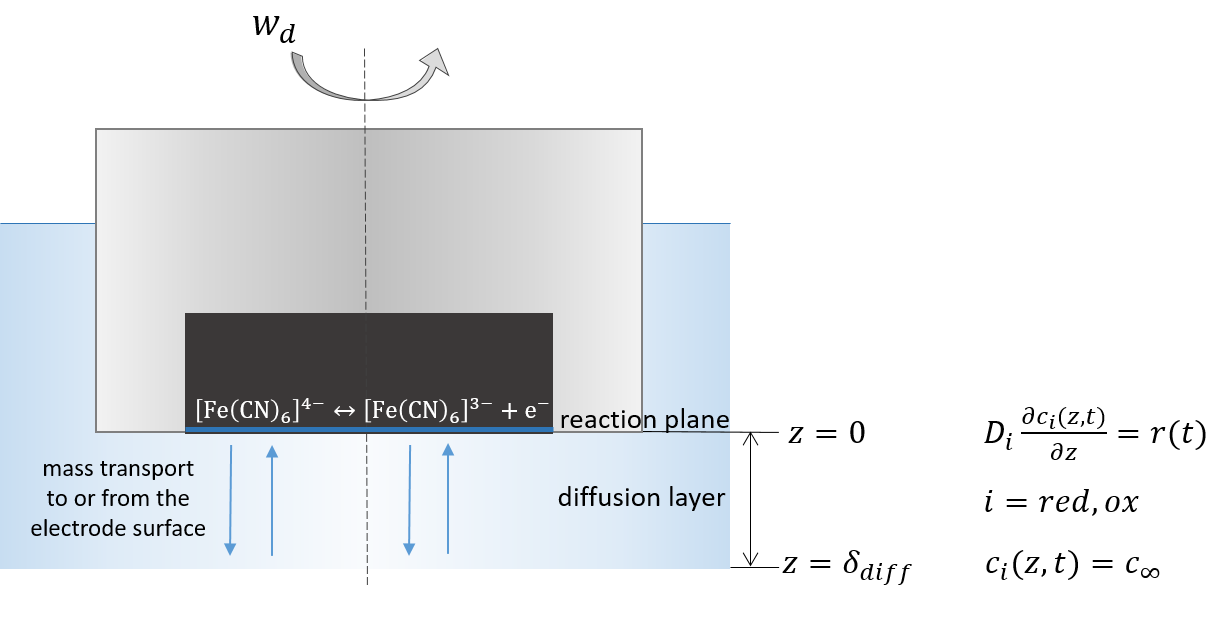}
	\caption{Schematic of coupled electrochemical kinetics and diffusion model}
	\label{fig:ferr_rotation}
\end{figure}

\label{subsubsec:battery}
This section presents the performance of the parametric DMD on a model of the ferrocyanide redox reaction. The reaction kinetics under the influence of the rotation rate of the rotating disc electrode is of interest \cite{VidPAetal11}. A schematic representation of the investigated system is shown in \Cref{fig:ferr_rotation}. This reaction can be considered as a model reaction with coupled electrochemical kinetics and mass transport. Similar to the first battery model, the governing equations of this model are based on mass and charge conservation laws as well.

The mass conservation law is described by the second Fick's law assuming that convective terms can be neglected, see \cref{eq:fc_fick}.

\begin{equation}
\label{eq:fc_fick}
\frac{\partial c_{i}(z, w_d, t)}{\partial t}=D_{i} \frac{\partial^{2} c_{i}(z,w_d, t)}{\partial z^{2}}, \quad i = red, ox, 
\end{equation}
where the subscript $i$ stands either for the reduced (ferrocyanide, \ch{[Fe(CN)6]\mch[4]}) or oxidized (ferricyanide, \ch{[Fe(CN)6]\mch[3]}) form, and $c_{i}, D_i$ are their corresponding concentration and diffusion coefficients, respectively.

The charge balance can be described as, 
\begin{equation}
\label{eq:fc_charge}
C_{dl} \frac{dE(w_d,t)}{dt} = J(w_d,t) - Fr(w_d,t),
\end{equation}
where $E(w_d,t)$ is the electrode potential, $C_{dl}$ is the double-layer capacitance, $J(w_d,t)$ is the cell current density, $F$ is Faraday constant,  and $r(w_d,t)$ is the nonlinear reaction rate, computed by Butler-Volmer kinetics,
\begin{equation}
\label{eq:ferr_reactionrate}
r(w_d,t)= k\left\{\frac{c_{\mathrm{red}}(0, w_d,t)}{c_{\mathrm{red}, \infty}} \exp \left(\beta f\left(E(w_d,t)-E_{\mathrm{r}}\right)\right)-\frac{c_{\mathrm{ox}}(0,w_d, t)}{c_{\mathrm{ox}, \infty}} \exp \left(-(1-\beta) f\left(E(w_d,t)-E_{\mathrm{r}}\right)\right)\right\}.
\end{equation}

Here, $E_r$ is the equilibrium electrode potential, $\beta$ is the charge transfer coefficient, and $f$ is determined as $F/RT$, where $T$ is the temperature, and $R$ is the universal gas constant.

The main source of the nonlinearity comes from $r(t)$ and its coupling with the diffusion of the reacting species through the boundary conditions (given in \Cref{fig:ferr_rotation} and \cref{eq:bc_ferr}). 
\begin{equation}
\label{eq:bc_ferr}
\left.D_i \frac{\partial c_{\mathrm{i}}(z,w_d, t)}{\partial z}\right|_{z=0}=\pm r(w_d,t), \quad i=\text { red or ox }.
\end{equation}
In \Cref{tab:parameter_ferr}, we list all the important parameters used to construct the model and their ranges of change.


\begin{table}
  \caption{Parameters in the model for ferrocyanide reaction.}
  \label{tab:parameter_ferr}
  \begin{center}
    \begin{tabular}{l||c|c}
      Parameters & Variables & Value Range\\
      \hline
      rotation rate & $w_d$ ($rpm$) & $[500, 5000]$ \\
      input potential & $E_{total}(V)$ & $[-0.4, 1.0]$\\
      double layer capacity & $C_{dl}$ ($F/m^2$) & $0.2$\\
      charge transfer coefficient & $\beta$ & $ 0.5 $\\
      reaction rate constant & k ($m/s$) &$[10^{-7}, 10^{-2}]$\\
      ohmic resistance of the electrolyte & $R_{\Omega}$($\Omega$) &$[1,100]$\\
      diffusivity coefficient for the ferrocyanide & $D_{red}$ ($m^2/s$) & $[10^{-10}, 9 \times 10^{-10}]$ \\
      diffusivity coefficient for the ferricyanide & $D_{ox}$ ($m^2/s$) & $[10^{-10}, 9 \times 10^{-10}]$\\
      \hline
    \end{tabular}
  \end{center}
\end{table}

\begin{figure}
\begin{subfigure}{0.5\textwidth}
\includegraphics[width=0.95\linewidth]{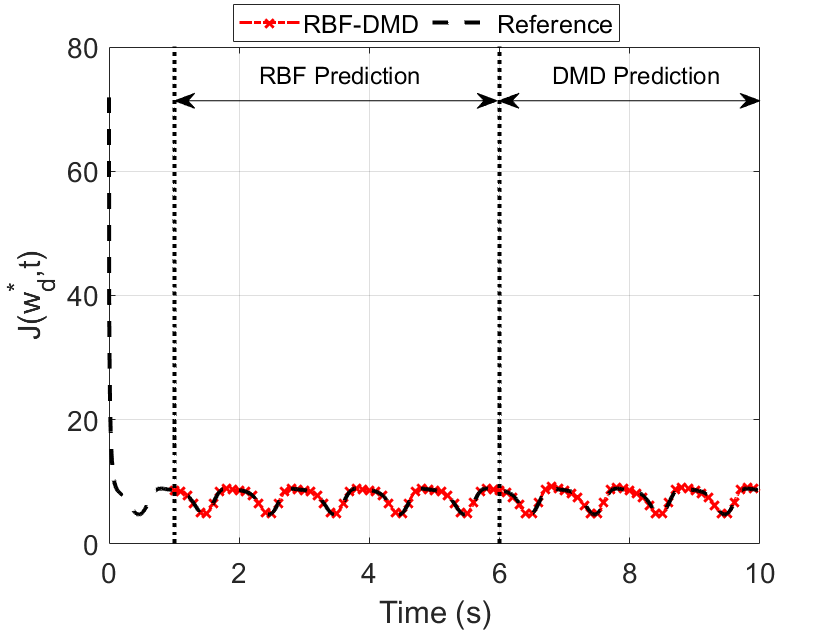} 
\caption{}
\label{fig:FERR_model_sol_1}
\end{subfigure}
\begin{subfigure}{0.5\textwidth}
\includegraphics[width=0.95\linewidth]{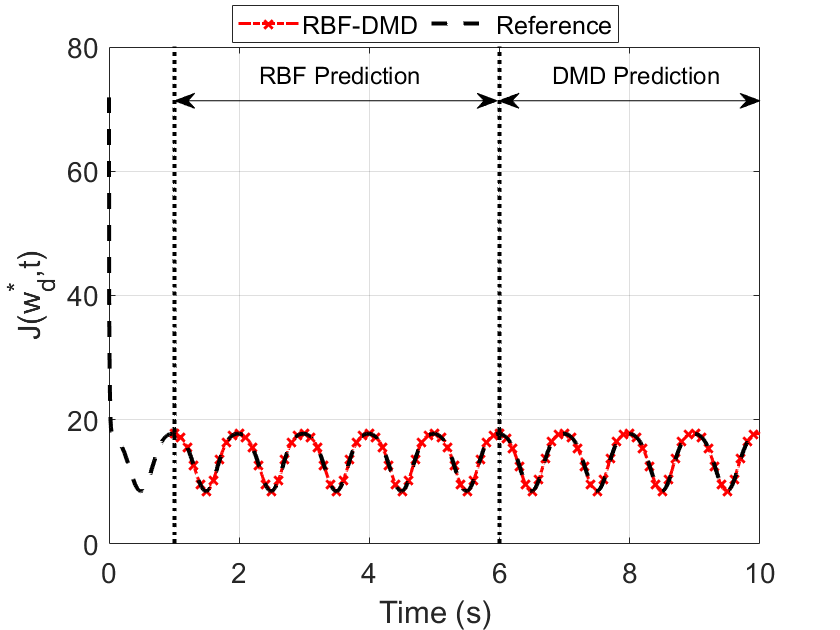}
\caption{}
\label{fig:FERR_model_sol_4}
\end{subfigure}
\caption{Ferrocyanide reaction model: the parametric DMD solution vs the reference solution for the current density $J(w_d^*,t)$. (a) $ w_d^* = 1000 \,rpm$. (b) $ w_d^*  = 4800 \,rpm$.}
\label{fig:FERR_model_sol}
\end{figure}

\begin{figure}
\begin{subfigure}{0.5\textwidth}
\includegraphics[width=0.95\linewidth]{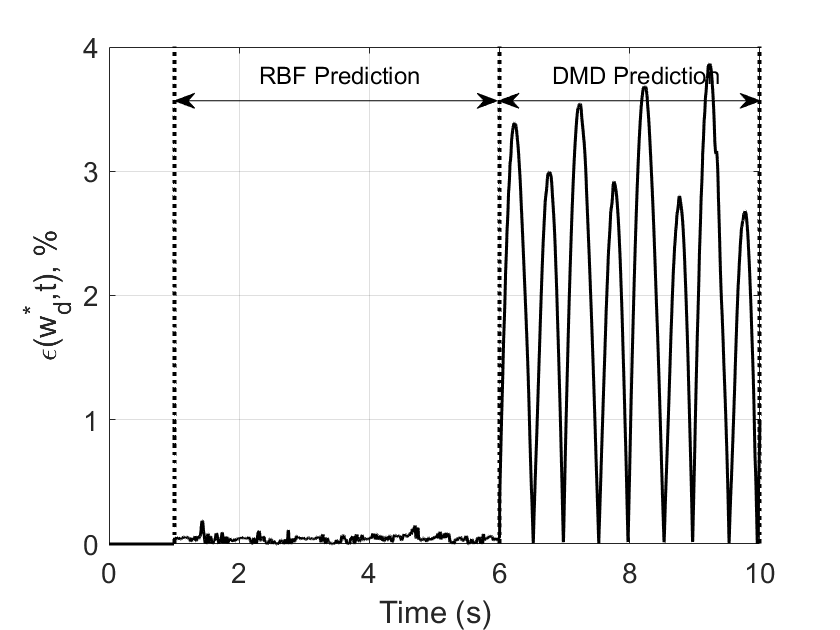} 
\caption{}
\label{fig:FERR_model_error_1}
\end{subfigure}
\begin{subfigure}{0.5\textwidth}
\includegraphics[width=0.95\linewidth]{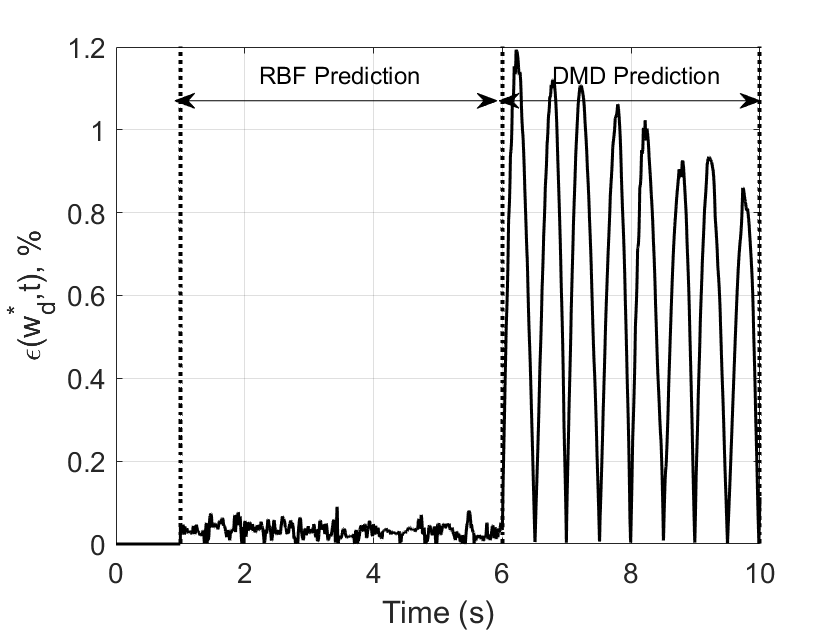}
\caption{}
\label{fig:FERR_model_error_4}
\end{subfigure}
\caption{Ferrocyanide reaction model: the relative error of the parametric DMD solution for the current density $J(w_d^*,t)$. (a) $ w_d^* = 1000 \,rpm$. (b) $ w_d^* = 4800 \,rpm$.}
\label{fig:FERR_model_error}
\end{figure}

\begin{figure}
\centering
\includegraphics[width=0.5\linewidth]{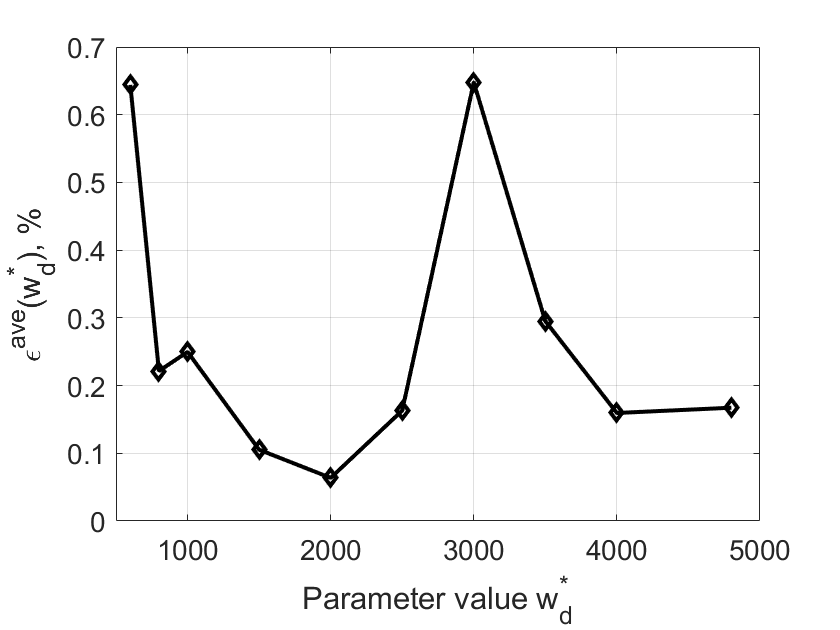}
\caption{Ferrocyanide reaction model: the time-average relative error of the parametric DMD solution for the current density $J(w_d^*,t)$ at different testing rotation rates $w_d^*$.}
\label{fig:FERR_model_merror}
\end{figure}

We study the influence of the rotation rate $w_d$ (\Cref{fig:ferr_rotation}) on the system output (current density $J(w_d,t)$). The rotation
speed $w_d$ of the rotating disc electrode determinates the thickness of diffusion layer $\delta_{diff}$ for ox or red, as shown below:
\begin{equation}
\label{eq:diffusionlayer}
\delta_{diff}=1.61 D_i^{1/3} \nu^{1/6} \w_d^{-1/2}, \quad i=\text { red or ox },
\end{equation}
where $\nu$ is the kinematic viscosity. The thickness of the diffusion layer further has impacts on the concentration in \cref{eq:fc_fick} and its boundary as $c_i\left(\delta_{diff}, t\right)=c_{i, \infty}, \quad i=$ red or ox (see also \Cref{fig:ferr_rotation}). The equations in \cref{eq:fc_charge} and \cref{eq:ferr_reactionrate} are discretized in space using finite differences. The dimension of the discretized system is $n = 4003$ while the simulation time is set as $10s$ with 10 periods. $20$ different rotation rates as training parameters are uniformly sampled in the range of $[500, 5000] \,rpm$. The kernel DMD with Gaussian kernel in our parametric DMD method (\Cref{algorithm:kerneldmd}) is selected in this example. $r$ in step 2 of \Cref{algorithm:kerneldmd} is chosen according to the criteria in \cref{eq:truncationtolerance} with $\eta = 0.5 \%$. \Cref{fig:FERR_model_sol} presents the current density computed by the parametric DMD and the reference solution. The relative error changing with time at two testing samples of $w_d$ and the time average relative errors at 10 different testing rotation rates are plotted in \Cref{fig:FERR_model_error} and \Cref{fig:FERR_model_merror}, respectively. In \Cref{fig:FERR_model_error}, the relative errors at all time instances are below $5 \%$. In \Cref{fig:FERR_model_merror}, the time average relative error at all testing $w_d^*$ samples, i.e., $w_d^* = 600, 800, 1000, 1500, 2000, 2500, 3000, 3500, 4000, 4800 \,rpm$ are under $0.7 \%$ when using parametric DMD. 


The runtime comparison for this model is listed in \Cref{tab:FERR_computation_time}. At the offline stage, generating snapshots and training RBF network take 1452.29 seconds and 0.082 seconds, respectively. The online runtime is the average value over 10 testing samples of different rotation rates. RBF predicts the current in $[0, T_0]$ using 0.082 seconds and in $[T_0 , T]$,  DMD uses 7.312 seconds. The total runtime at the online stage is around 7.4 seconds, which is much less than that of solving the original system (FOM simulation) by an ODE solver with 77.703 seconds.

\begin{table}[h]
{\footnotesize
  \caption{Ferrocyanide reaction model: The computation time (seconds) of parametric DMD and that of the FOM simulation.} \label{tab:FERR_computation_time}
\begin{center}
  \begin{tabular}{|c|c|c|c|c|} \hline
  Snapshot generation & RBF training & RBF prediction & DMD  prediction & FOM simulation\\ \hline
  1452.290 &  2.453  & 0.719 & 6.202 & 99.825 \\ \hline
\end{tabular}
\end{center}}
\end{table}

\subsection{FitzHugh–Nagumo model}%
\label{subsubsec:FHNmodel}
We further consider the nonlinear Fitz–Hugh Nagumo model as a benchmark example used in many existing works \cite{morAsiABetal21, morBenB15, morChaS10, morBenGG18}. This model is designed to simulate the spike generation in an excitable system, for example in a neuron. The describing equations read:

\begin{equation}
\begin{gathered}
\varepsilon v_{t}(x, \varepsilon, t)=\varepsilon^{2} v_{x x}(x, \varepsilon, t)+f(v(x, \varepsilon, t))-w(x, \varepsilon, t)+c \\
w_{t}(x,\varepsilon,  t)=b v(x, \varepsilon, t)-\gamma w(x, \varepsilon, t)+c \\
y(x,\varepsilon, t) = [v(0,\varepsilon,t), w(0,\varepsilon,t)]^T,
\end{gathered}
\end{equation}
with $f(v) = v(v-0.1)(1-v)$ as the cubic nonlinear term and the boundary conditions are:
\begin{equation*}
\begin{array}{llr}
v(x, \varepsilon,0)=0, & w(x,\varepsilon, 0)=0, & x \in[0,L], \\
v_x(0,\varepsilon, t)=-i_o(t), & v_x(L, \varepsilon,t)=0, & t \geq 0,
\end{array}
\end{equation*}

The unknown state variable, $v(x,\varepsilon,t)$ is the membrane potential, and $w(x,\varepsilon,t)$ is a recovery of the potential. Parameters are $ b,c, \varepsilon$ and $\gamma$. In this numerical test, the operating parameter is $\varepsilon$, changing from $0.02$ to $0.03$, while other parameters are fixed as $L = 20$, $b = 0.5$, $c= 0.05$ and $\gamma = 2$. The input term is $i_o(t) = 50000t^3e^{-15t}$. The output vector $y(x,\varepsilon,t)\in \mathbb{R}^{2}$ includes two outputs: the membrane potential and the recovery of the potential at the left boundary.

After discretization by the finite difference method, the resulting ODE is solved by the ODE solver ode15s in MATLAB. The total number of states is $n = 16384$. The time span is $[0, 10]s$ with the time step $\delta t = 0.01s$. The snapshots are taken in the time interval $[0,8]s$. The number of the equidistant samples in $[0.02, 0.03]$ in the training phase is $15$. For this example, kernel DMD is chosen in \Cref{algorithm:kerneldmd}.

\begin{figure}
\begin{subfigure}{0.5\textwidth}
\includegraphics[width=0.99\linewidth]{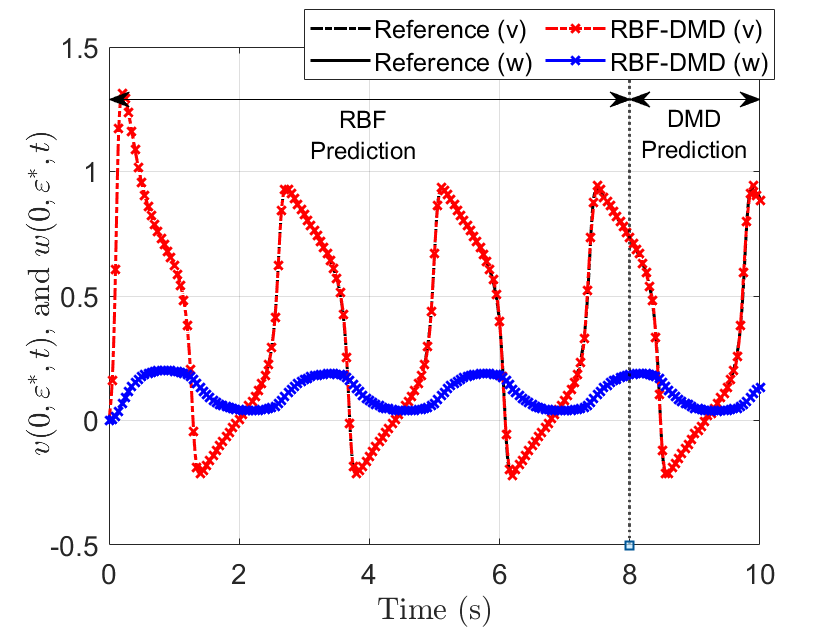} 
\caption{}
\label{fig:FHN_model_sol_0225_1}
\end{subfigure}
\begin{subfigure}{0.5\textwidth}
\includegraphics[width=0.99\linewidth]{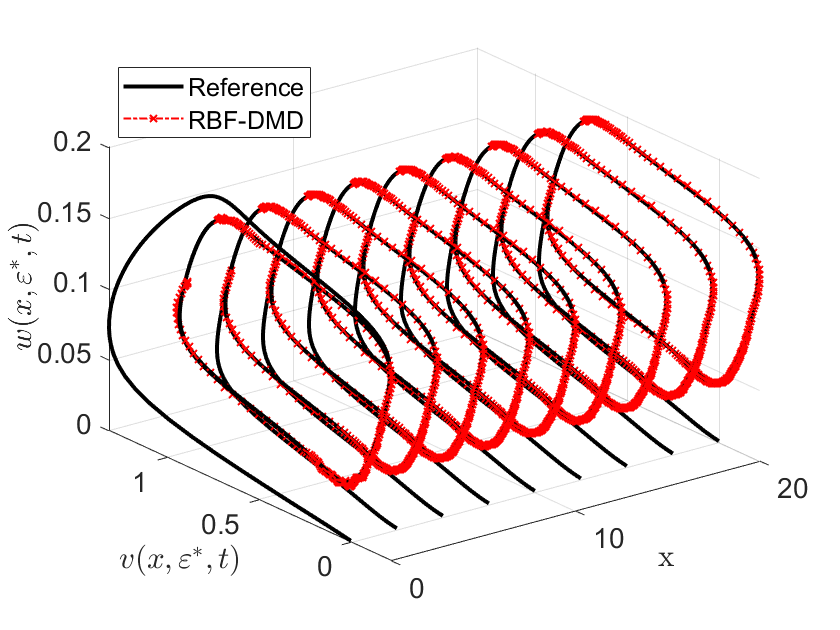}
\caption{}
\label{fig:FHN_model_sol_0225_2}
\end{subfigure}
\begin{subfigure}{0.5\textwidth}
\includegraphics[width=0.99\linewidth]{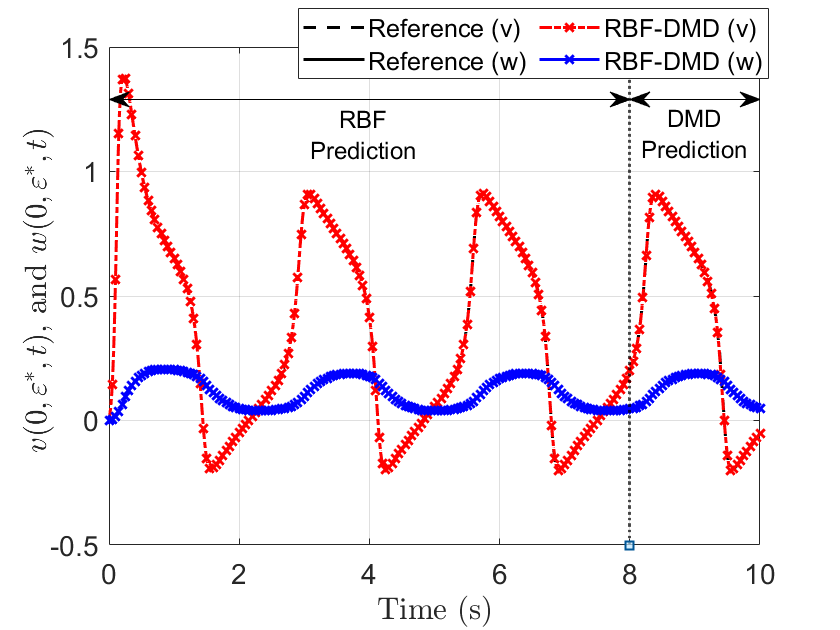} 
\caption{}
\label{fig:FHN_model_sol_0275_1}
\end{subfigure}
\begin{subfigure}{0.5\textwidth}
\includegraphics[width=0.99\linewidth]{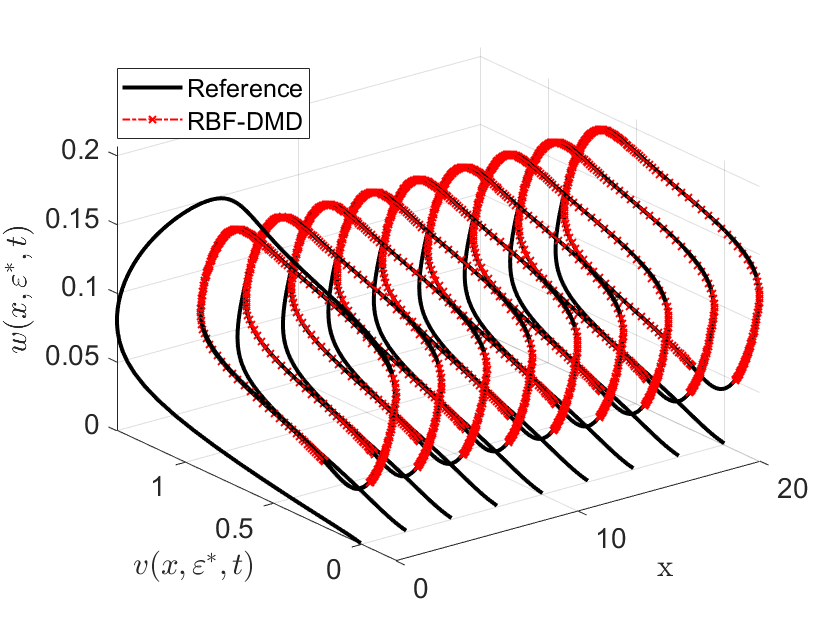}
\caption{}
\label{fig:FHN_model_sol_0275_2}
\end{subfigure}
\caption{FitzHugh-Nagumo model: the parametric DMD solution vs the reference solution. (a) The outputs $v(0,\varepsilon^*,t)$ and $w(0,\varepsilon^*,t)$ when $\varepsilon^* = 0.0225$. (b) Limit cycles of $v(x,\varepsilon^*,t)$ and $w(x,\varepsilon^*,t)$ when $\varepsilon^* = 0.0225$. (c) The outputs $v(0,\varepsilon^*,t)$ and $w(0,\varepsilon^*,t)$ when $\varepsilon^* = 0.0275$. (d) Limit cycles of $v(x,\varepsilon^*,t)$ and $w(x,\varepsilon^*,t)$ when $\varepsilon^* = 0.0275$.}
\label{fig:FHN_model_sol}
\end{figure}

\begin{figure}
\begin{subfigure}{0.5\textwidth}
\includegraphics[width=0.99\linewidth]{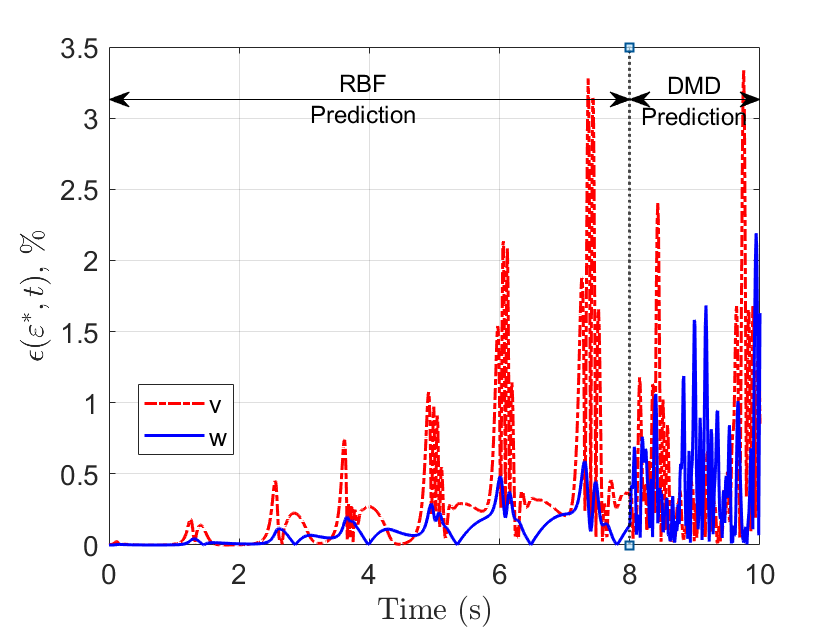} 
\caption{}
\label{fig:FHN_model_sol_3}
\end{subfigure}
\begin{subfigure}{0.5\textwidth}
\includegraphics[width=0.99\linewidth]{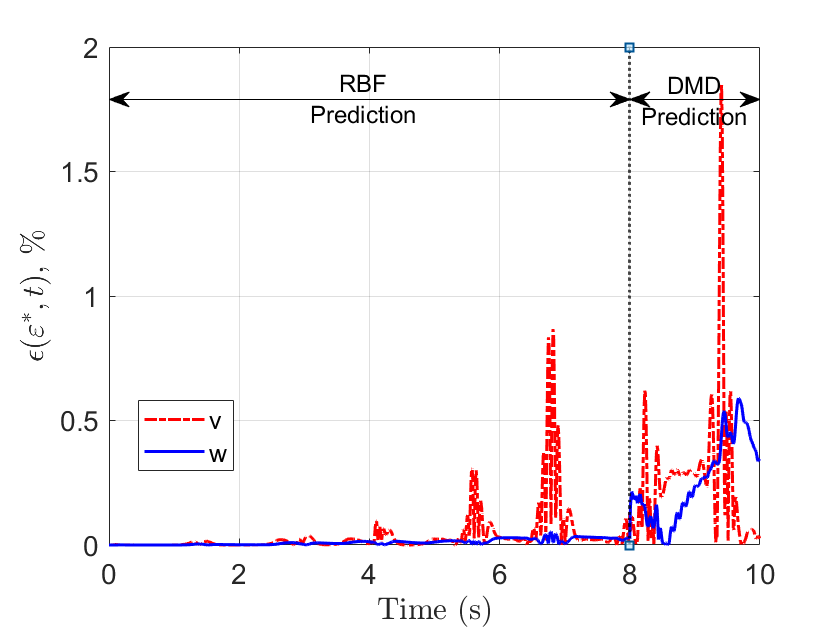}
\caption{}
\label{fig:FHN_model_sol_4}
\end{subfigure}
\caption{FitzHugh-Nagumo model: the relative error of the parametric DMD for $v(0,\varepsilon^*,t)$ and $w(0,\varepsilon^*,t)$. (a) $\varepsilon^* = 0.0225$. (b) $\varepsilon^* = 0.0275$.}
\label{fig:FHN_model_err}
\end{figure}

\begin{figure}
\centering
\includegraphics[width=0.5\linewidth]{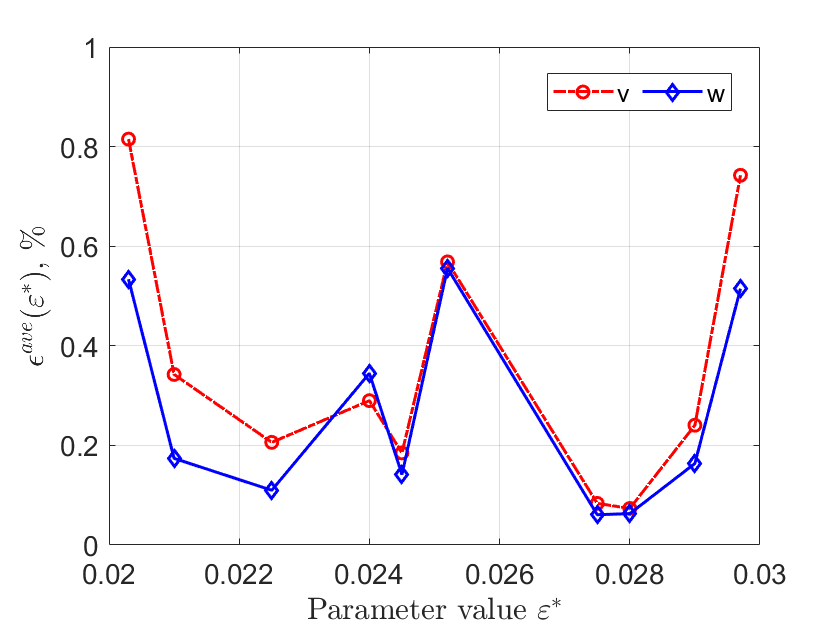}
\caption{FitzHugh-Nagumo model: the time-average relative error of the parametric DMD solution at different testing $\varepsilon^*$.}
\label{fig:FHN_model_merr}
\end{figure}

The numerical results are shown in \Cref{fig:FHN_model_sol}. \Cref{fig:FHN_model_sol_0225_1} and \Cref{fig:FHN_model_sol_0275_1} show the evolution of the two outputs $v(0,\varepsilon^*,t)$ and $w(0,\varepsilon^*,t)$ when $\varepsilon^* = 0.0225$ and $\varepsilon^* = 0.0275$. As is shown in these figures, at the online stage of the proposed parametric DMD, RBF first predicts the solution at the testing $\varepsilon^*$ in the time interval $[0, 8]s$, then DMD predicts the evolution of the solution in the time period $[8, 10]s$. The red line is the parametric DMD results for $v(0,\varepsilon^*,t)$ and the blue line stands for $w(0,\varepsilon^*, t)$. Both lines fit quite well with the black reference solution. \Cref{fig:FHN_model_sol_0225_2} and \Cref{fig:FHN_model_sol_0275_2} are their corresponding phase-space diagrams. \Cref{fig:FHN_model_err} is the relative error changing with time when $\varepsilon^* = 0.0225$ and $\varepsilon^* = 0.0275$. The maximum relative error of these two cases is around $2-3.5  \% $.  \Cref{fig:FHN_model_merr} is the time average of the relative errors over all testing parameters, i.e., $\varepsilon^* = 0.021, 0.0225, 0.024, 0.0245, 0.0252,0.027, 0.0275, 0.029 $. Their values never exceeds $1 \%$ in all these testing cases. Through these error plots, it can be confirmed that the proposed method works well for this nonlinear dynamic system.

The computation time is also listed in \Cref{tab:FHN_computation_time}. 12656.380 seconds are needed at the offline stage for generating training snapshots for training samples. RBF training takes 1.395 seconds based on these training samples. At the online stage, the RBF prediction costs 0.5435 seconds and the DMD prediction costs 0.694 seconds. Solving the original full order model for a single testing parameter takes around 1057.4 seconds. It can be observed that around 850 times speed-up is achieved at the online phase when using parametric DMD.

\begin{table}[h]
{\footnotesize
  \caption{FitzHugh-Nagumo model: The computation time (seconds) of parametric DMD and that of the FOM simulation.}  \label{tab:FHN_computation_time}
\begin{center}
  \begin{tabular}{|c|c|c|c|c|} \hline
   Snapshot generation & RBF training & RBF prediction & Online prediction & FOM simulation\\ \hline

 12656.380 & 1.395 & 0.5435  & 0.694 & 1057.400 \\ \hline
\end{tabular}
\end{center}}
\end{table}

\section{Conclusion}
\label{sec:conclu}
We propose a non-intrusive parametric model order reduction method combining the DMD and RBF. When heavy computations are needed for multi-query tasks in the parametric case, especially for predicting the nonlinear dynamics, the proposed parametric DMD is promising for prediction in both the parameter and the time domain.

The proposed method is tested on several examples and their results are compared with the reference solutions obtained by direct simulations of the original models. The results demonstrate that the proposed algorithm is effective. For the P2D battery model, where the frequency of the current is the changing parameter, parametric DMD predicts the output potential at a new frequency with high accuracy. The second example of the ferrocyanide redox reaction is parametrized with rotation rates. The numerical results also indicate the high accuracy of the parametric DMD. The FitzHugh-Nagumo model further manifests the effective reduction and acceptable accuracy of the parametric DMD for the large nonlinear dynamic system.

Further improvements can be done in several directions. Firstly, all the numerical examples are based on a single parameter and the data are from simulation. Parametric DMD could also be applied to real experimental data with multiple parameters, which is of high interest in the design of experiments (DoE). Secondly, DMD and its related topic are being developed with a rapid speed, the proposed method could be further extended to new variants of DMD.

\section*{Data and code availability}
\addcontentsline{toc}{section}{Data and code availability}
Data and code will be available in a public repository later, and on request.

\section*{Declarations of interest}
\addcontentsline{toc}{section}{Declarations of interest}
None.

\section*{Funding}
\addcontentsline{toc}{section}{Funding}
This research is partially supported by the International Max Planck Research School for Advanced Methods in Process and Systems Engineering (IMPRS ProEng), Magdeburg, Germany.





\bibliographystyle{elsarticle-num} 
\bibliography{refs_new}

\end{document}